\documentstyle[12pt]{article}
\input {cyracc.def}
\font\ququ=cmr10 scaled \magstep1
\font\tencyr=wncyr7 scaled \magstep1
\def\rus{\tencyr\cyracc}

\newcommand{\reb}[1]{\mbox{\bf  (\ref{#1})}}
\newcommand{\re}[1]{\mbox{\rm  (\ref{#1})}}

\newenvironment{proof}
{\noindent {\sl Proof.}\quad }{\hfill
$\square$ \vskip1.1ex\noindent }

\newenvironment{proof*}
{\noindent {\sl Proof.}\quad }{\hfill
$\square$}

\renewcommand{\theequation}{\thesection .\arabic{equation}}

\renewcommand{\thesubsubsection}{\theequation .\arabic{subsubsection}}

\catcode`\@=\active
\catcode`\@=11
\def\@eqnnum{\hbox to .01pt{}\rlap{\bf \hskip -\displaywidth(\theequation)}}
\catcode`\@=12

\newenvironment{s}[1]
{ \vskip1.2ex \refstepcounter{equation}
\noindent {\bf \theequation\quad #1.} \begin{sl}}{\end{sl}
\vskip1.1ex\noindent }

\newenvironment{subs}[1]
{\vskip1.2ex \refstepcounter{equation}
\noindent {\bf (\theequation)\quad #1.} }{\quad}

\newenvironment{subsubs}[1]
{\vskip1.2ex \refstepcounter{subsubsection}
\noindent {\bf (\thesubsubsection)\quad #1.} }{\quad}

\newcommand {\sekt}[1]
{{\vskip2.5ex\refstepcounter{section}\setcounter{equation}{0}
\noindent\large \bf \thesection\quad \parbox[t]{424pt}{#1}
\nopagebreak\vskip1.5ex\noindent}}

\newcommand {\g}{{\frak g}}

\newcommand {\me}{{\frak m}}
\newcommand {\n}{{\frak n}}
\newcommand {\p}{{\frak p}}

\newcommand {\te}{{\frak t}}
\newcommand {\tri}{{\frak sl}_2}

\newcommand {\z}{{\frak z}}
\newcommand {\esi}{\varepsilon}
\newcommand {\ap}{\alpha}

\newcommand {\vp}{\varphi}

\newcommand {\V}{{\Bbb V}}
\newcommand {\hb}{{\bf h}}

\newcommand {\eb}{{\bf{e}}}

\newcommand {\mathrm}{\rm\mbox}

\newcommand {\ad}{{\mathrm{ad\,}}}
\newcommand {\Ad}{{\mathrm{Ad\,}}}

\newcommand {\Ker}{{\mathrm{Ker\,}}}
\newcommand {\Ima}{{\mathrm{Im\,}}}
\newcommand {\rk}{{\mathrm{rk\,}}}

\newcommand {\GR}[2]{{\mathrm{{\bf #1}}}_{#2}}

\newcommand {\un}{\underline}

\newcommand {\vno}[1]{\vskip#1 ex\noindent}
\newcommand {\rar}{\rightarrow}
\newcommand {\qus}{\hfill $\square$ \vno{1.1}}

\newcommand {\beq}{\begin{equation}}
\newcommand {\eeq}{\end{equation}}

\newfam\Bbbfam\newfam\eufam%
\font\Bbbfont=msbm10 scaled 1200%
\font\olala=msam10 scaled 1200%
\font\frak=eufm10 scaled 1400%
\font\Bbbsmallfont=msbm7 scaled 1000%
\textfont\Bbbfam=\Bbbfont\scriptfont\Bbbfam=\Bbbsmallfont
\font\euzw=eufm10 scaled 1200%
\font\euac=eufm7 scaled 1200%
\textfont\eufam=\euzw\scriptfont\eufam=\euac%
\def\frak{\fam\eufam}%
\def\Bbb{\fam\Bbbfam}%

\def\varnothing{\hbox {\Bbbfont\char'077}}
\def\square{\hbox {\olala\char"03}}
\def\Bbbk{\hbox {\Bbbfont\char'174}}

\begin{document}
\setlength{\parskip}{2pt plus 4pt minus 0pt}
\hfill {\scriptsize September 10, 1999}
\vskip1ex

\noindent
{\Large \bf A classification of the principal nilpotent pairs
\vno{1.3}in simple Lie algebras and related problems}
\bigskip \\
{\bf Alexander G. Elashvili\footnote{supported
in part by Grant INTAS-OPEN-97-1570} and
Dmitri I. Panyushev\footnote{supported in part by
R.F.F.I. Grant {\rus N0}\,98--01--00598} }

\medskip
\smallskip

\noindent{\large \bf Introduction}
\vno{2}%
Let $\g$ be a semisimple Lie algebra over an algebraically closed field $\Bbbk$
of characteristic zero and $G$ its adjoint group.
The notion of a {\it principal nilpotent pair\/}
is a double
counterpart of the notion of a regular (= principal) nilpotent element in
$\g$. Roughly speaking, a principal nilpotent pair $\eb=(e_1,e_2)$ consists
of two commuting elements in $\g$ that can independently be contracted to
the origin and such that their simultaneous centralizer
has the minimal possible dimension, i.e., $\rk\g$.
The very definition and the basic results are due to V.\,Ginzburg
\cite{vitya}. He showed that the theory of principal nilpotent pairs
yields a refinement of well-known results due to B.\,Kostant
on regular nilpotent
elements in $\g$ and has interesting applications to Representation Theory.
In particular, he proved that the number of $G$-orbits of principal nilpotent
pairs is finite and gave a classification for $\g={\frak sl}(\V)$.
Trying to achieve a greater generality, Ginzburg also introduced a wider class of
{\it distinguished nilpotent pairs}  (see precise definitions in sect.\,1)
and, again, classified them for ${\frak sl}(\V)$. Although finiteness for
${\frak sl}(\V)$ follows from the classification, it is not  known in general
whether the number of $G$-orbits of distinguished nilpotent pairs
is finite.
\\[.6ex]
Our aim is to present a classification of
the distinguished and principal nilpotent pairs in the classical simple Lie
algebras. Because a classification of {\sl principal\/} nilpotent pairs in the
exceptional simple Lie algebras was obtained in \cite{wir}, we thus
have a complete classification of such pairs in all simple algebras.
As a by-product,
we obtain finiteness for the number of orbits of the distinguished
pairs in all classical Lie algebras. However the finiteness problem
remains open for the exceptional simple algebras.
\\[.5ex]
Roughly speaking, the reason for success is that a classical Lie algebra
has the tautological representation. Let $\V$ be the space of this
representation.
The classification is given in terms of combinatorial objects which are
called {\it skew-graphs}. A skew-graph $\Gamma$
is a finite oriented graph whose set
of nodes is a subset of ${\Bbb Q}\times\Bbb Q$, see precise definition in
\reb{skew}. To attach a skew-graph to a distinguished pair $\eb$, we exploit
the bi-grading of $\V$ determined by a {\it characteristic\/} $\hb$ of $\eb$
(see sect.\,1 for characteristics). This is an extension of Ginzburg's
approach to ${\frak sl}(\V)$.
For each classical series, we describe the appropriate class of skew-graphs and
construct a mapping from the set of $G$-orbits of distinguished pairs to
this class of skew-graphs. This mapping is a bijection unless
$\g={\frak so}_{2n}$ and the target skew-graph $\Gamma$ is connected. In the
last case, the fibre over $\Gamma$ consists of two $SO_{2n}$-orbits.
Everybody who has heard about
`very even' partitions for $\g={\frak so}_{2n}$ (see e.g. \cite[5.1]{CoMc})
will not find it surprising. For principal pairs, the classification is of
the same form, only appropriate classes of skew-graphs become smaller.
The classification is summarized in Theorem \ref{MAIN}.
\\[.5ex]
The notion of a connected skew-graph is close to that of a
{\it skew-diagram\/} considered by Ginzburg in connection with distinguished
pairs in ${\frak sl}(\V)$. But our connected skew-graphs are more
formalized objects and
we found it more convenient to use the language of graphs in the situation,
where the combinatorial object can be disconnected.
An advantage of using coordinates for the nodes of $\Gamma$ is two-fold:
1) one distinguishes easily the classes of skew-graphs attached to
different classical series and 2) one can immediately read off a formula for
a characteristic of a distinguished pair in question.
\\[.5ex]
In section 1, we recall the necessary definitions and results. In section 2,
we introduce the ``appropriate" classes of skew-graphs and describe
the corresponding
distinguished and principal pairs. The classification itself is contained in
section 3. We show how to associate a skew-graph to a distinguished pair
and that the resulting skew-graphs are exactly those considered in section 2.
In section 4, we give an explicit description of the centralizer for
all principal pairs and describe {\sl rectangular\/} distinguished pairs.
\\[.5ex]
{\bf Notation.} The ground field $\Bbbk$ is algebraically closed and of
characteristic zero. Algebraic groups are denoted by capital Roman letters,
whereas their Lie algebras by the corresponding small Gothic letters.
Throughout, 
$\g$ is a semisimple Lie algebra and $G$ is its adjoint group.
For any set $M\subset\g$, we let $\z_\g(M)$ (resp. $Z_G(M)$) denote
the centralizer of $M$ in $\g$ (resp. in $G$).
For $M=\{a,\dots,z\}$, we simply write $\z_g(a,\dots,z)$ or $Z_G(a,\dots,z)$.
If $x\in\g$ and $s\in G$, we write $s{\cdot}x$
in place of $(\Ad s)x$. \par
If $\V$ is a $\g$-module and $M\subset\g$, then $\V^M$ denotes the
subspace of $M$-fixed
vectors. If $\V$ is equipped with a non-degenerate  bilinear form $(~,~)$
and $W\subset \V$, then $W^\perp=\{v\in\V\mid (v,w)=0\quad \forall w\in W\}$.
\par
$\langle a,\dots,z\rangle$ is the linear span of
elements of a vector space; \par
${\Bbb P}=\{0,1,2,\dots\}$, ${\Bbb N}=\{1,2,\dots\}$.

\sekt{Preliminaries on nilpotent pairs\nopagebreak}%
We begin with recalling some definitions from
\cite{vitya}.
\par
$\bullet$
A pair $\eb=(e_1,e_2) \in \g\times\g$
is called {\it nilpotent\/} in $\g$, if \\
(i) $[e_1, e_2]=0$  and \quad
(ii) for any $(t_1,t_2)\in \Bbbk^*\times \Bbbk^*$,  there
exists $g=g(t_1,t_2)\in G$  such that
$(\,t_1e_1,\,t_2e_2\,)=(g{\cdot}e_1,\,g{\cdot}e_2\,).$
\par
$\bullet$
A pair of semisimple elements $\hb=(h_1,h_2)\in\g\times\g$ is called
an {\it associated semisimple pair\/} for $\eb$, if $h_1,h_2$ have rational
eigenvalues in $\g$ and
\beq [h_1, h_2]=0,\quad [h_i, e_j]=\delta_{ij} e_j \quad
(i,j\in \{1,2\})\ .  \label{comrel}
\eeq
\indent
$\bullet$
A nilpotent pair is called {\it principal\/},
if $\dim\z_\g(\eb)=\rk\g$.
\par
$\bullet$
A nilpotent pair is called {\it pre-distinguished\/}, if
$\z_\g(\eb)$ contains no semisimple elements.
\par
$\bullet$
A nilpotent pair is called {\it distinguished\/}, if it is pre-distinguished
and there exists an associated semisimple pair $\hb$
such that $\z_\g(\hb)$ is a Cartan subalgebra.
\\[.5ex]
It was shown in \cite{vitya} that each nilpotent pair has an associated
semisimple pair.
The main tool in studying nilpotent pairs is the bi-grading
$\g=\bigoplus_{p,q\in {\Bbb Q}}\g_{p,q}$,
where $\g_{p,q}=\{x\in\g\mid [h_1,x]=px,[h_2,x]=qx\}$.
The pairs $(p,q)$ with $\g_{p,q}\ne 0$
will be referred to as the {\it eigenvalues\/} of $\hb$ in $\g$. An
eigenvalue $(p,q)$ is said to be {\it integral}, if $p,q\in\Bbb Z$. Otherwise
it is called {\it fractional}. The same terminology is used for
the corresponding eigenspaces. Of course, this can also be performed for
any  $\g$-module in place of $\g$. Dealing with distinguished and principal
nilpotent pairs we shall usually omit the adjective `nilpotent'.\\
The following result was proved in \cite{vitya}.
\begin{s}{Theorem {\ququ (Ginzburg)}} \label{princ} \\
Let $\hb$ be an associated semisimple pair for a principal pair
$\eb$. Then \par
{\sf (i)} $\z_\g(\hb)$ is a Cartan subalgebra; \par
{\sf (ii)} the eigenvalues of $\hb$ in $\g$ are integral; \par
{\sf (iii)} $\displaystyle\z_\g(\eb)=\bigoplus_{p,q\in{\Bbb P}, (p,q)\ne (0,0)}
\z_\g(\eb)_{p,q}$.
\end{s}%
It follows from the theorem that a principal pair is distinguished.
To utilize the bi-grading determined by $\hb$ for obtaining internal properties
of a nilpotent pair, one
need to have a kind of conjugacy theorem for all associated semisimple pairs.
In \cite{vitya}, such a theorem was proved for the
pre-distinguished pairs. However, to get the result in full
generality, the definition of an associated semisimple pair was to be modified.
This was done in \cite{IMRN}.
\par
$\bullet$
A pair $\hb$ of semisimple elements is called
a {\it characteristic\/} of $\eb$, if it satisfies \re{comrel} and
$h_1,h_2\in\z_\g(\eb)^\perp$.
\\[.5ex]
Note that for the pre-distinguished pairs the second condition is automatically
satisfied. It was shown in \cite[1.4]{IMRN} that any nilpotent
pair has a characteristic, which is unique within to conjugacy, and $h_1,
h_2$ are rational semisimple elements. The following assertion will
repeatedly be used in course of our classification of distinguished pairs.
We repeat the proof given in [loc.\,cit, 1.5].
\begin{s}{Lemma} \label{distin}
Let $\eb$ be a nilpotent pair in $\g$. Suppose \par
{\sf (i)} $\hb$ is a semisimple pair satisfying Eq.\,\re{comrel}; \par
{\sf (ii)} $\z_\g(\hb)\cap\z_\g(\eb)=0$.\\
Then $\eb$ is pre-distinguished.
\end{s}\begin{proof}
Let $\g=\z_\g(\hb)\oplus\me$ be the $Z_G(\hb)$-stable decomposition. Then
$\z_\g(\eb)\subset\me$ and therefore $\z_\g(\eb)$ is orthogonal to
$h_1,h_2$. That is, $\hb$ is a characteristic of $\eb$. Consider the
algebraic Lie algebra $\n_\g(\eb)=
\{x\in\g\mid [x,e_i]\in\langle e_i\rangle\quad i=1,2\}$. Then
$\n_\g(\eb)=\z_\g(\eb)\oplus\langle h_1,h_2\rangle$ as $\Bbbk$-vector space.
Using algebraic Levi decompositions for $\n_\g(\eb)$ and $\z_\g(\eb)$,
one easily shows that $\z_\g(\hb)$ contains $\z_\g(\eb)^{red}$,
a reductive Levi subalgebra of $\z_\g(\eb)$ (see \cite[1.4]{IMRN}).
Hence $\z_\g(\eb)^{red}=0$ and we are done.
\end{proof}

\sekt{The examples of distinguished and principal pairs in classical Lie
algebras\nopagebreak}%
In this section we describe examples of distinguished and principal pairs in
the classical simple Lie algebras. It will be shown in section~3 that
these examples actually exhaust all such pairs.
\begin{subs}{Skew-graphs and skew-diagrams} \label{skew} \end{subs}
\\
A {\it  skew-graph\/} $\Gamma$ is a finite oriented graph
satisfying the following conditions: \par
{\sf (i)} the set of nodes
$n(\Gamma)$ is a subset of ${\Bbb Q}\times{\Bbb Q}$:
\par
{\sf (ii)} the barycentre of the set $n(\Gamma)$ is in the origin,
i.e., $\sum_{(i,j)\in n(\Gamma)}i=0$ and $\sum_{(i,j)\in n(\Gamma)}j=0$;\par
{\sf (iii)} the set of arrows $a(\Gamma)$ contains elements
of only two types: 
$\bullet\!\rar\!\bullet$
\quad or \quad
\begin{tabular}{@{}c@{}}\vbox{\hbox{$\bullet$\rule{0ex}{1ex}}
\hbox{$\uparrow$\rule{0ex}{2ex}\strut}
\hbox{$\bullet$\rule{0ex}{1ex}}}
\end{tabular}\ , where the length of each arrow is 1;
\par
{\sf (iv)} if the  nodes $(i,j), (i+1,j+1)$ lie in the same connected component
$\Gamma_1$ of $\Gamma$, then $(i,j+1), (i+1,j)$ also belong to $\Gamma_1$
and four possible arrows inside this
square belong to $a(\Gamma_1)$.
\\[.6ex]
Different components are allowed to have
some common nodes.
\\[.6ex]
It is easily seen that, in a connected skew-graph, the set $n(\Gamma)$
entirely determines $a(\Gamma)$. For this reason, we shall sometimes identify
$\Gamma$ and $n(\Gamma)$ in such a case and depict a connected skew-graph
$\Gamma$ as collection of squares of size 1 on the plane, an element
$(i,j)\in n(\Gamma)$ being
the centre of the corresponding square. The corresponding object is said to be
a {\it skew-diagram} $\Gamma$.
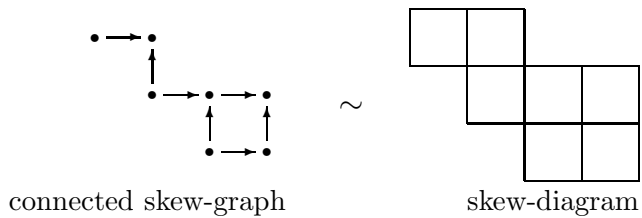
\begin{center}
\begin{picture}(52,40)(0,15) \setlength{\unitlength}{0.03in}
\multiput(0,20)(10,-10){3}{\circle*{1.5}}
\multiput(10,20)(10,-10){3}{\circle*{1.5}}
\put(30,10){\circle*{1.5}}
\multiput(2,20)(10,-10){3}{\vector(1,0){6}}
\multiput(10,12)(10,-10){2}{\vector(0,1){6}}
\put(30,2){\vector(0,1){6}}
\put(22,10){\vector(1,0){6}}
\put(-15,-10){{\small connected skew-graph}}
\end{picture}
\qquad\quad $\sim$ \qquad
\begin{picture}(45,40)(0,15) \setlength{\unitlength}{0.03in}
\put(-5,25){\line(1,0){20}}
\put(-5,15){\line(1,0){40}}
\put(5,5){\line(1,0){30}}
\put(15,-5){\line(1,0){20}}
\put(-5,15){\line(0,1){10}}
\put(5,5){\line(0,1){20}}
\put(15,-5){\line(0,1){30}}
\put(25,-5){\line(0,1){20}}
\put(35,-5){\line(0,1){20}}
\put(5,-10){{\small skew-diagram}}
\end{picture}
\end{center}
\vskip6ex\noindent
Associated with these two types of objects, there are two different languages
that will freely be used in the sequel. Thus, we gain an advantage to use those
language that is more suitable in a specific situation.
Usually, it is not hard to switch
from one language to another.
For instance, a southwest (resp. northeast) corner of a skew-diagram $\Gamma$
corresponds to a source (resp. sink) of the oriented graph $\Gamma$.
A skew-diagram is said to be a {\bf sw}- (resp. {\bf ne}-) {\it Young
diagram\/}, if
it has a unique southwest (resp. northeast) corner. The {\bf sw}- and
{\bf ne}-Young diagrams together are called Young diagrams. The corresponding
skew-graphs are said to be {\it Young graphs}.
\\[.5ex]
The following statement was proved in \cite{vitya}. We included a proof for
two reasons: (a) convenience of the reader and (b) it is the prototype of
our argument for the other classical series.
\begin{s}{Proposition {\ququ (series {\bf A})}} \label{existA}
\\
{\sf 1}. Any connected skew-graph $\Gamma$ with $n$ nodes determines a
distinguished pair in ${\frak sl}_n$. The $SL_n$-orbit of this pair is uniquely
determined by $\Gamma$. \\
{\sf 2.} Any Young graph with $n$ nodes determines a principal
pair in ${\frak sl}_n$.
\end{s}\begin{proof}
1. Identify the nodes of $\Gamma$ with the elements of a
basis of an $n$-dimensional $\Bbbk$-vector space $\V$.
If $(i,j)\in n(\Gamma)$, then $v_{i,j}$ stands
for the corresponding basis vector in $\V$. It is convenient to assume that
$v_{i,j}=0$, if $(i,j)\not\in n(\Gamma)$.
Define $e_1$ (resp. $e_2$)
to be the operator corresponding to the horizontal (resp. vertical) arrows
in $a(\Gamma)$: If $(i,j)\in n(\Gamma)$ then
$e_1(v_{i,j}):=v_{i+1,j}$ and $e_2(v_{i,j}):=v_{i,j+1}$.
Obviously, $e_1,e_2$ are nilpotent endomorphisms of $\V$.
>From \ref{skew}(iv) it follows that $[e_1,e_2]=0$.
Define the semisimple endomorphisms of $\V$ by the formulae:
$h_1(v_{i,j})=iv_{i,j}$, $h_2(v_{i,j})=jv_{i,j}$.
The condition on the barycentre of $n(\Gamma)$ is equivalent to that
$\mbox{tr\,}(h_i)=0$, $i=1,2$.
It is easily seen that $h_1,h_2,e_1,e_2$ satisfy
relations \re{comrel}. Since $\Gamma$ is connected, a straightforward
computation shows that $\z_{{\frak sl}(\V)}(\hb)$ is the set of
all traceless endomorphisms of $\V$ that are diagonal in the basis
$\{v_{i,j}\}_{(i,j)\in n(\Gamma)}$. It is then easily seen that
$\z_{{\frak sl}(\V)}(\hb)\cap\z_{{\frak sl}(\V)}(\eb)=0$.
Therefore $\eb$ is distinguished by Lemma~\ref{distin}.
\par
To prove uniqueness, assume that $\{w_{i,j}\}_{(i,j)\in n(\Gamma)}$ is
another basis of $\V$ and $\eb'=(e'_1,e'_2)$ is the corresponding
distinguished pair. Clearly, there exist $s\in SL(\V)$ such that
$sv_{i,j}=\ap w_{i,j}$ for all $(i,j)$ and some $\ap\in\Bbbk^*$.
Then $s{\cdot}e_i=e'_i$ ($i=1,2$).
\par
2. Let $\Gamma$ be a {\bf sw}-Young graph and $(i_0,j_0)$ the source of it.
Then $v_{i_0,j_0}$ is a cyclic vector in $\V$ relative to $e_1,e_2$.
Therefore $\z_{{\frak sl}(\V)}(\eb)$ is generated by all non-zero powers
$e_1^ke_2^l$, $(k,l)\ne (0,0)$. It is easily seen that
$e_1^ke_2^l\ne 0$ if and only if $(i_0+k,j_0+l)\in n(\Gamma)$.
Whence $\dim\z_{{\frak sl}(\V)}(\eb)=\dim\V-1=
\rk{\frak sl}(\V)$. \par
If $\Gamma$ is an {\bf ne}-Young graph, we may consider the dual
endomorphisms $e_1^*,e_2^*$ of $\V^*$. Obviously, the pair
$\eb^*=(e_1^*,e_2^*)$ corresponds to the {\bf sw}-Young diagram $\Gamma^t$ and
$\dim\z_{{\frak sl}(\V)}(\eb)=\dim\z_{{\frak sl}(\V^*)}(\eb^*)$.
\end{proof}%
\begin{subs} {Centrally-symmetric skew-graphs} \label{cs}
\end{subs}
Let us say that a skew-graph $\Gamma$ is {\it centrally-symmetric} (=\,c.-s.),
if the sets $n(\Gamma)$ and $a(\Gamma)$ are centrally-symmetric in the usual
sense (the latter being considered without orientation).
For a connected $\Gamma$, this implies that the coordinates of all nodes
are in $\frac{1}{2}{\Bbb Z}\times\frac{1}{2}{\Bbb Z}$.
All connected c.-s. skew-graphs fall into three classes. Namely, $\Gamma$ is
said to be
\vskip1ex
\halign{\indent#&\quad#\hfil&\quad#\cr
{\it integral\/}, &if $i,j\in\Bbb Z$ &for a (any) $(i,j)\in n(\Gamma)$;\cr
{\it semi-integral\/}, &if $i+j\in\frac{1}{2}+{\Bbb Z}$
&for a (any) $(i,j)\in n(\Gamma)$;\cr
{\it non-integral\/}, &if $i,j\in\frac{1}{2}+{\Bbb Z}$
&for a (any) $(i,j)\in n(\Gamma)$.\cr}
\vskip1ex\noindent
A geometric property distinguishing these classes is the position of centre
of symmetry with respect to nodes and arrows.
As we shall see, these classes naturally arise in describing distinguished
pairs in ${\frak so}(\V)$ and ${\frak sp}(\V)$.
Notice that there are two different {\it sorts\/} of semi-integral
skew-graphs: one
with $i\in\Bbb Z$ and another with $j\in\Bbb Z$. Clearly, $\# n(\Gamma)$ is
odd if $\Gamma$ is integral and is even for the other two classes.
We say that a connected skew-graph
is {\it rectangular\/}, if the corresponding skew-diagram is a rectangle.
Obviously, a rectangular skew-graph is c.-s.
\begin{s}{Proposition {\ququ (series {\bf B})}} \label{existB} \\
{\sf 1.} Let $\Gamma$ be a c.-s. skew-graph with at most two connected
components $\Gamma_0$, $\Gamma_1$ such that $\Gamma_0$ is integral and
$\Gamma_1$ is non-integral (it is allowed that $\Gamma_1=\varnothing$).
Then such $\Gamma$ determines a distinguished
pair in ${\frak so}(\V)$, where $\dim\V=\# n(\Gamma)$ is odd.
The $SO(\V)$-orbit of this pair is uniquely determined by $\Gamma$;
\\
{\sf 2.} Any rectangular integral skew-graph $\Gamma$ determines a principal
nilpotent pair in ${\frak so}(\V)$, where $\dim\V=\# n(\Gamma)$ is odd.
\end{s}{\sl Proof.\quad}
We follow the same way as in \re{existA}, with
necessary alterations. \\
1. Obviously, each $\Gamma_i$ is c.-s., $n(\Gamma_0)\cap n(\Gamma_1)=
\varnothing$, and $\# n(\Gamma)$ is odd.
Identify the nodes of $\Gamma_i$ with the elements of a basis of a
$\Bbbk$-vector space $\V_i$ ($i=0,1$). Set $\V=\V_0\oplus\V_1$,
$\g={\frak so}(\V)$, and $\g_i={\frak so}(\V_i)$.
Define the {\sl symmetric\/} bilinear
form on $\V$ by $(v_{i,j},v_{-k,-l}):=\delta_{ik}\delta_{jl}$.
As above, $e_1$ (resp. $e_2$) will correspond to the
horizontal (resp. vertical) arrows of $\Gamma$. However, to obtain
$e_1, e_2$ respecting the bilinear form, one must include signs in their
definitions. The problem is to choose a distribution of signs such that
$e_1, e_2$ respect the bilinear form and still commute. A solution is as
follows. For $(i,j)\in n(\Gamma)$, set
\beq  \label{so-basis}
e_1(v_{i,j})=(-1)^{i+j}v_{i+1,j} \quad \mbox{and}
\quad e_2(v_{i,j})=(-1)^{i+j+1}v_{i,j+1}\ .
\eeq
(Note that $i{+}j$ is integral for both components of $\Gamma$.)
Then  $e_1,e_2\in{\frak so}(\V)$ and $[e_1,e_2]=0$.
Define $h_1,h_2$ by the same formulae as in the proof of
Prop.~\ref{existA}. It is immediate that $h_1,h_2\in {\frak so}(\V)$
and relations \re{comrel} are satisfied.
There is an $\hb$-equivariant isomorphism
${\frak so}(\V)\simeq\wedge^2(\V_0\oplus\V_1)\simeq
{\frak so}(\V_0)\oplus{\frak so}(\V_1)\oplus (\V_0\otimes\V_1)$.
Since each $h_i$ has integral eigenvalues in $\V_0$ and non-integral ones in
$\V_1$, we obtain $(\V_0\otimes\V_1)^{\langle h_1,h_2\rangle}=0$.
Therefore $\z_\g(\hb)=
\z_{\g_0}(\hb)\oplus\z_{\g_1}(\hb)$. Since $\Gamma_i$ is connected,
$\z_{\g_i}(\hb)$ is a Cartan subalgebra of $\g_i$ ($i=0,1$).
Hence $\z_\g(\hb)$ is a Cartan subalgebra of $\g$. It consists of all operators
in $\g$ that are diagonal in the basis $\{v_{i,j}\}_{(i,j)\in n(\Gamma)}$.
It then easy to verify that
$\z_{\g}(\hb)\cap\z_{\g}(\eb)=0$.
Thus, $\eb$ is distinguished by Lemma~\ref{distin}. \par
To prove uniqueness, assume that $\{w_{i,j}\}_{(i,j)\in n(\Gamma)}$ is
another basis of $\V$, satisfying the condition $(w_{i,j},w_{-k,-l})=
\delta_{ik}\delta_{jl}$, and $\eb'=(e'_1,e'_2)$ is the corresponding
distinguished pair. There exist $s\in O(\V)$ such that
$sv_{i,j}=w_{i,j}$ for all $(i,j)$. Then $s{\cdot}e_i=e'_i$ ($i=1,2$).
If $\mbox{det\,}s=1$, we are done. If $\mbox{det\,}s=-1$, then one makes
a correcting alteration: Set $s'v_{i,j}=-w_{i,j}$ for all
$(i,j)\in n(\Gamma_0)$ and $s'v_{i,j}=w_{i,j}$ for all $(i,j)\in n(\Gamma_1)$.
Then $s'{\cdot}e_i=e'_i$ and since $\# n(\Gamma_0)$ is odd,
$\mbox{det\,}s'=1$.
\par
2. A rectangular skew-graph has a unique source, say $(i_0,j_0)$, which
provides us with a cyclic vector in $\V$ relative to $e_1,e_2$. Therefore
$\z_{{\frak so}(\V)}(\eb)$ is generated by all non-zero powers
$e_1^ke_2^l$ lying in ${\frak so}(\V)$. The last condition is equivalent to
that $k+l$ is odd. Since $\Gamma$ is rectangular integral c.-s., we have \\
\hbox to 450pt {\hfil $\displaystyle
\#\{(k,l)\mid k+l\ \mbox{ is odd }\ \&\ (i_0+k,j_0+l)\in n(\Gamma)\}=
\frac{\# n(\Gamma)-1}{2}=\rk{\frak so}(\V)$ . \hfil $\square$}
\\[1ex]
{\bf Remark.} Notice that this proposition applies in degenerate cases as well.
For instance, if $\Gamma$ is a horizontal chain, then it has no vertical
arrows and therefore $e_2=0$. Since $e_1$ has a single block in Jordan normal
form, it is a regular nilpotent element in ${\frak so}(\V)$.
This provides us with a trivial principal pair $(e_1,0)$.
\begin{s}{Proposition {\ququ (series {\bf C})}} \label{existC} \\
{\sf 1.}  Let $\Gamma$ be a c.-s. skew-graph with at most two connected
components $\Gamma_0$, $\Gamma_1$ that are both
semi-integral, of different sorts (it is allowed that $\Gamma_1=\varnothing$).
Then such $\Gamma$ determines a distinguished
pair in ${\frak sp}(\V)$, where $\dim\V=\# n(\Gamma)$ is even.
The $Sp(\V)$-orbit of this pair is uniquely determined by $\Gamma$;
\\
{\sf 2.} A rectangular semi-integral skew-graph $\Gamma$ determines a
principal
nilpotent pair in ${\frak sp}(\V)$, where $\dim\V=\# n(\Gamma)$ is even.
\end{s}{\sl Proof.\quad}
1. Since $\Gamma_0$ and $\Gamma_1$ are semi-integral and of different sorts,
$n(\Gamma_0)\cap n(\Gamma_1)=\varnothing$. For definiteness sake, we assume
that $i\in\Bbb Z$ (resp. $j\in\Bbb Z$) if $(i,j)\in n(\Gamma_0)$
(resp. $(i,j)\in n(\Gamma_1)$).
Identify the nodes of $\Gamma_i$ with the elements of a basis of a
$\Bbbk$-vector space $\V_i$ ($i=0,1$). As above, $v_{i,j}$ stands for
the basis vector corresponding to $(i,j)\in n(\Gamma)$.
Set $\V=\V_0\oplus\V_1$,
$\g={\frak sp}(\V)$, and $\g_i={\frak sp}(\V_i)$.
Define the {\sl alternating\/} bilinear
form on $\V$ by
\begin{center}
$(v_{i,j},v_{-k,-l}):=0$ unless $i=k$ and $j=l$; \\
if $(i,j)\in n(\Gamma_0)$, then $(v_{i,j},v_{-i,-j}):=1$ for $j>0$;\\
if $(i,j)\in n(\Gamma_1)$, then $(v_{i,j},v_{-i,-j}):=1$ for $i>0$.
\end{center}
Formulae for $e_1,e_2$ also become more involved.
For  $(i,j)\in n(\Gamma_0)$, we set:
\beq  \label{sp-basis} \qquad
e_1(v_{i,j})=\left\{\begin{array}{rl} -v_{i+1,j}, & \mbox{if } j>0 \\
v_{i+1,j}, & \mbox{if } j<0 \end{array} \right.\ ,\qquad
 e_2(v_{i,j})=\left\{\begin{array}{rl}
v_{i,j+1}, & \mbox{if } j>0 \\
(-1)^iv_{i,j+1},  & \mbox{if } j=-1/2 \\
-v_{i,j+1}, & \mbox{if } j+1<0
\end{array} \right. .
\eeq
The formulas for $(i,j)\in n(\Gamma_1)$ are obtained in a `transposed' fashion:
\[  
\qquad
e_1(v_{i,j})=\left\{\begin{array}{rl}
v_{i+1,j}, & \mbox{if } i>0 \\
(-1)^jv_{i+1,j},  & \mbox{if } i=-1/2 \\
-v_{i+1,j}, & \mbox{if } i+1<0
\end{array} \right.\ , \qquad
e_2(v_{i,j})=\left\{\begin{array}{rl} -v_{i,j+1}, & \mbox{if } i>0 \\
v_{i,j+1}, & \mbox{if } i<0 \end{array} \right. .
\]
It is not hard to verify that $[e_1,e_2]=0$ and
$e_1,e_2 \in {\frak sp}(\V)$.
Define $h_1,h_2$ by the same formulae as in the proof of
Prop.~\ref{existA}. It is immediate that $h_1,h_2\in {\frak sp}(\V)$
and relations \re{comrel} are satisfied.
There is an $\hb$-equivariant isomorphism
${\frak sp}(\V)\simeq{\cal S}^2(\V_0\oplus\V_1)\simeq
{\frak sp}(\V_0)\oplus{\frak sp}(\V_1)\oplus (\V_0\otimes\V_1)$.
Since $\Gamma_0$ and $\Gamma_1$ are of different sorts,
we obtain $(\V_0\otimes\V_1)^{\langle h_1,h_2\rangle}=0$.
Therefore $\z_\g(\hb)=
\z_{\g_0}(\hb)\oplus\z_{\g_1}(\hb)$. Since $\Gamma_i$ is connected,
$\z_{\g_i}(\hb)$ is a Cartan subalgebra of $\g_i$ ($i=0,1$).
Hence $\z_\g(\hb)$ is a Cartan subalgebra of $\g$.
It then easy to verify that
$\z_{\g}(\hb)\cap\z_{\g}(\eb)=0$.
Thus, $\eb$ is distinguished by Lemma~\ref{distin}.
\par
The proof of uniqueness is the same as in Prop.~\ref{existB} and even simpler,
because the full isometry group of an alternating form is connected.
\par
2. This part of the proof is almost the same as in \ref{existB}(2). The only
distinction is in the last formula. Since $\Gamma$ is {\sl semi\/}-integral,
we have \\
\hbox to 450pt {\hfil $\displaystyle
\#\{(k,l)\mid k+l\ \mbox{ is odd }\ \&\ (i_0+k,j_0+l)\in n(\Gamma)\}=
\frac{\# n(\Gamma)}{2}=\rk{\frak sp}(\V)$ . \hfil $\square$}
\\[1ex]
As the description in the even-dimensional orthogonal case is longer,
we split it in two propositions: one about distinguished pairs and another
about principal pairs. A new phenomenon appearing here is that different
connected components of $\Gamma$ may have a common node.
\begin{s}{Proposition {\ququ (series {\bf D})}} \label{existDd} \\
{\sf (i)} Any connected non-integral c.-s. skew-graph $\Gamma$  determines a
distinguished pair in ${\frak so}(\V)$, where $\dim\V=\# n(\Gamma)$ is even; \\
{\sf (ii)} Let $\Gamma$ be a skew-graph with two connected components
$\Gamma_1$ and $\Gamma_2$ that are c.-s. and integral. Then $\Gamma$ determines
a distinguished pair in ${\frak so}(\V)$ in the following two cases:
\par
a) \ $n(\Gamma_2)=\{(0,0)\}$; \par
b) \ $\# n(\Gamma_i)>1$ $(i=1,2$) and
$n(\Gamma_1)\cap n(\Gamma_2)=\{(0,0)\}$.
\\
{\sf (iii)} Let $\Gamma$ be a skew-graph with three connected components
$\Gamma_i$ ($i=0,1,2$) such that $\Gamma_0$ is as in (i) and
$\Gamma_1$, $\Gamma_2$ are as in (ii).  Then $\Gamma$ represents a
distinguished pair in $\V$, where $\dim\V=\sum_i \# n(\Gamma_i)$; \\
{\sf (iv)} The $SO(\V)$-orbit of distinguished pair is uniquely determined by
$\Gamma$ if and only if $\Gamma$ is disconnected. In case $\Gamma$ is
connected (see (i)), two different orbits arise.
\end{s}\begin{proof}
(i) The argument here is analogous to that in \re{existB}.
Consider a basis $\{v_{i,j}\}$ of an even-dimensional $\Bbbk$-vector space
$\V$ indexed by the nodes of $\Gamma$.
The symmetric bilinear form on $\V$ is defined by
$(v_{i,j},v_{-k,-l}):=0$ unless $i=k$ and $j=l$,
and $(v_{i,j},v_{-i,-j}):=1$.
For any $(i,j)\in n(\Gamma)$, set
\[
e_1(v_{i,j})=(-1)^{i+j}v_{i+1,j}\ , \quad
e_2(v_{i,j})=(-1)^{i+j+1}v_{i,j+1}
\]
\[  \mbox{and}\quad
h_1(v_{i,j})=iv_{i,j},\quad h_2(v_{i,j})=jv_{i,j}\ .
\]
The verification of all required relations is easy. \par
(ii) Let $\V_i$ be an odd-dimensional $\Bbbk$-vector space with a basis
indexed by the nodes of $\Gamma_i$, $i=1,2$.
According to Prop.~\ref{existB}(1), each
$\Gamma_i$ represents a distinguished pair $(e_1^{(i)},e_2^{(i)})$,
with characteristic $(h_1^{(i)},h_2^{(i)})$, in ${\frak so}(\V_i)$.
Then straightforward computations show that the direct sum
$(e_1,e_2):=(e_1^{(1)}+e_1^{(2)},e_2^{(1)}+e_2^{(2)})$ is a
distinguished pair in $\g={\frak so}(\V_1\oplus\V_2)={\frak so}(\V)$ in both
cases:
\par
(a) Here $\dim\V_2=1$ and $e_1^{(2)}=e_2^{(2)}=h_1^{(2)}=h_2^{(2)}=
0$. Hence $e_i=e_i^{(1)}$ and $h_i=h_i^{(1)}$. Since
${\frak so}(\V)\simeq\wedge^2(\V_1\oplus\V_2)\simeq\wedge^2\V_1\oplus
(\V_1\otimes\V_2)\simeq {\frak so}(\V_1)\oplus\V_1$ as $h_1$- and $h_2$-module
and $\dim\V_1^{\langle h_1,h_2\rangle}=1$, we obtain
$\dim\z_{\g}(\hb)=\dim\z_{{\frak so}(\V_1)}(\hb)+1=
\rk{\frak so}(\V_1)+1=\rk{\frak so}(\V)$, i.e., $\z_{\g}(\hb)$
is a Cartan subalgebra. It is also easily seen that $\z_{\g}(\hb)
\cap\z_{\g}(\eb)=0$.
\par
(b) Here ${\frak so}(\V)\simeq{\frak so}(\V_1)\oplus{\frak so}(\V_2)\oplus
(\V_1\otimes\V_2)$. Hence $\z_\g(\hb)=
{\frak so}(\V_1)^{\langle h_1,h_2\rangle}\oplus
{\frak so}(\V_2)^{\langle h_1,h_2\rangle}\oplus
(\V_1\otimes\V_2)^{\langle h_1,h_2\rangle}$. The assumption
$n(\Gamma_1)\cap n(\Gamma_2)=\{(0,0)\}$ implies that
$(\V_1\otimes\V_2)^{\langle h_1,h_2\rangle}=
\V_1^{\langle h_1,h_2\rangle}\otimes\V_2^{\langle h_1,h_2\rangle}$
is 1-dimensional.
Therefore $\dim\z_{\g}(\hb)=\rk{\frak so}(\V_1)+
\rk{\frak so}(\V_2)+1=\rk\g$, i.e., $\z_\g(\hb)$ is a Cartan subalgebra.
For $\z_\g(\eb)$, we have the similar formula
\[
\z_\g(\eb)=
{\frak so}(\V_1)^{\langle e_1,e_2\rangle}\oplus
{\frak so}(\V_2)^{\langle e_1,e_2\rangle}\oplus
(\V_1\otimes\V_2)^{\langle e_1,e_2\rangle} \ .
\]
It follows from the connectedness of $\Gamma_i$ that
${\frak so}(\V_i)^{\langle e_1,e_2\rangle}\cap
{\frak so}(\V_i)^{\langle h_1,h_2\rangle}=0$ ($i=1,2$).
Finally, as both skew-graphs have more than 1 node, the space
$(\V_1\otimes\V_2)^{\langle h_1,h_2\rangle}$ is not $\eb$-stable.
Whence $\z_\g(\eb)\cap\z_\g(\hb)=0$.
We conclude by Lemma~\ref{distin}.
\par
(iii) Let $\V_i$ ($i=0,1,2$) be a $\Bbbk$-vector space with a basis
indexed by the nodes of $\Gamma_i$ and $\V:=\oplus_i\V_i$. The spaces
$\{\V_i\}$ are supposed to be pairwise orthogonal. Construct distinguished
pairs in ${\frak so}(\V_0)$ and ${\frak so}(\V_1\oplus\V_2)$
according to the recipes in (i) and (ii) respectively. Define $\eb$
to be the direct sum of these two pairs.
Again, we use the decomposition $\g=
{\frak so}(\V)\simeq{\frak so}(\V_0)\oplus{\frak so}(\V_1\oplus\V_2)\oplus
(\V_0\otimes(\V_1\oplus\V_2))$ to show that $\z_\g(\hb)=
\z_{{\frak so}(\V_0)}(\hb)\oplus\z_{{\frak so}(\V_1\oplus\V_2)}(\hb)$ and
hence $\z_\g(\hb)$ is a Cartan subalgebra of $\g$. The crucial equality here
is $(\V_0\otimes(\V_1\oplus\V_2))^{\langle h_1,h_2\rangle}=0$
(cf. proof of \ref{existB}(1)). Then the
problem of proving that $\z_\g(\hb)\cap\z_\g(\eb)=0$ becomes ``local", i.e.,
it reduces to $\V_0$ and $\V_1\oplus\V_2$, i.e., to preceding parts (i), (ii).
We conclude by Lemma \ref{distin}.
\par
(iv) As in Prop.~\ref{existB}, one proves that the distinguished pairs
constructed in (i),\,(ii), and (iii) form a single $O(\V)$-orbit. If
one of the connected components of $\Gamma$ has an odd number nodes, then
this orbit is actually an $SO(\V)$-orbit (cf. \ref{existB}). But if
$\Gamma$ is connected and hence contains an even number if nodes, then
this orbit is a union of two $SO(\V)$-orbits (the conjugating
transformation $s$ with $\mbox{det\,}s=-1$ can not be corrected).
\end{proof}%
Unlike the symplectic and odd-dimensional orthogonal case, principal
pairs in $\GR{D}{n}$ are also represented by some non-rectangular skew-graphs.
Before stating the result, we need to introduce this class of graphs.
\begin{subs}{Near-rectangular graphs and diagrams} \label{near}
\end{subs}
Let $\Gamma$ be a non-integral rectangular diagram. In this context,
non-integrality means that the length of each side of the rectangle is even.
Take the leftmost column of $\Gamma$ and remove from it either the lowest
(first) square or all the squares except for the highest (last) one.
Accordingly, execute the centrally-symmetric operation with the rightmost
column. This is one type of admissible transformations. Otherwise, the
similar thing can be done with the
first and last row of $\Gamma$. In any case, the diagram obtained
must be connected. The resulting
skew-diagrams (skew-graphs) are said to be {\it near-rectangular\/}.
Note that if $\Gamma$ is near-rectangular, then $\# n(\Gamma)\equiv 2
\pmod{4}$.
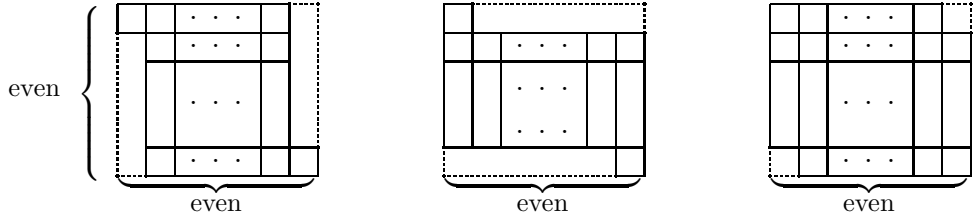
\begin{figure}[ht]
\caption{Near-rectangular diagrams}  \label{nears}
\setlength{\unitlength}{0.015in}
\centerline
{
\raisebox{36\unitlength}%
{$\left\{\parbox{1pt}{\vspace{60\unitlength}}\right.$}
\begin{picture}(80,77)
\put(-38,38){{\footnotesize even}}
\put(10,10){\line(1,0){60}}
\put(10,20){\line(1,0){60}}
\put(10,50){\line(1,0){50}}
\put(0,60){\line(1,0){60}}
\put(0,70){\line(1,0){60}}
\put(10,10){\line(0,1){60}}
\put(20,10){\line(0,1){60}}
\put(50,10){\line(0,1){60}}
\put(60,10){\line(0,1){60}}
\put(70,10){\line(0,1){10}}
\put(0,60){\line(0,1){10}}
\multiput(25,15)(0,20){3}{. . .}
\put(25,65){. . .}
\put(0,10){\dashbox{1}(0,50){}}
\put(0,10){\dashbox{1}(10,0){}}
\put(70,20){\dashbox{1}(0,50){}}
\put(60,70){\dashbox{1}(10,0){}}
\put(-.5,10){$\underbrace%
{\mbox{\hspace{68\unitlength}}}_{\mbox{{\footnotesize even}}}$}
\end{picture}
{\hspace{.4in}}
\begin{picture}(80,77)
\put(60,10){\line(1,0){10}}
\put(0,20){\line(1,0){70}}
\put(0,50){\line(1,0){70}}
\put(0,60){\line(1,0){70}}
\put(0,70){\line(1,0){10}}
\put(0,20){\line(0,1){50}}
\put(10,20){\line(0,1){50}}
\put(20,20){\line(0,1){40}}
\put(50,20){\line(0,1){40}}
\put(60,10){\line(0,1){50}}
\put(70,10){\line(0,1){50}}
\multiput(25,25)(0,15){3}{. . .}
\put(0,10){\dashbox{1}(60,0){}}
\put(0,10){\dashbox{1}(0,10){}}
\put(70,60){\dashbox{1}(0,10){}}
\put(10,70){\dashbox{1}(60,0){}}
\put(-.5,10){$\underbrace%
{\mbox{\hspace{68\unitlength}}}_{\mbox{{\footnotesize even}}}$}
\end{picture}
{\hspace{.4in}}
\begin{picture}(80,77)
\put(0,10){\dashbox{1}(0,10){}}
\put(0,10){\dashbox{1}(10,0){}}
\put(70,60){\dashbox{1}(0,10){}}
\put(60,70){\dashbox{1}(10,0){}}
\put(10,10){\line(1,0){60}}
\put(0,20){\line(1,0){70}}
\put(0,50){\line(1,0){70}}
\put(0,60){\line(1,0){70}}
\put(0,70){\line(1,0){60}}
\put(10,10){\line(0,1){60}}
\put(20,10){\line(0,1){60}}
\put(50,10){\line(0,1){60}}
\put(60,10){\line(0,1){60}}
\put(70,10){\line(0,1){50}}
\put(0,20){\line(0,1){50}}
\multiput(25,15)(0,20){3}{. . .}
\put(25,65){. . .}
\put(-.5,10){$\underbrace%
{\mbox{\hspace{68\unitlength}}}_{\mbox{{\footnotesize even}}}$}
\end{picture}
}
\end{figure}
\begin{s}{Proposition {\ququ (series {\bf D})}} \label{existDp} \\
{\sf (i)} Any non-integral rectangular or connected near-rectangular
skew-graph $\Gamma$ determines a principal
pair in ${\frak so}(\V)$, where $\dim\V=\# n(\Gamma)$ is even; \\
{\sf (ii)} Let $\Gamma$ be a skew-graph with two connected components
$\Gamma_1$ and $\Gamma_2$ that are c.-s. and integral. Then $\Gamma$ determines
a principal  pair in ${\frak so}(\V)$ in the following two cases: \par
a) $\Gamma_1$ is rectangular and $\# n(\Gamma_2)=1$; \par
b) $\Gamma_1$ is a horizontal chain and $\Gamma_2$ is a vertical chain.
\end{s}\begin{proof}
(i) For a non-integral rectangle, the proof is similar to that given in
Prop.~\ref{existB}(2), \ref{existC}(2) and is easy. Therefore we omit it. \\[1ex]
Suppose $\Gamma$ is a near-rectangular diagram. Let $\V$, the symmetric form,
and $\eb$ be the same as in \ref{existDd}(i).
The only thing that has to be verified is that
$\dim\z_{{\frak so}(\V)}(\eb)=\rk{\frak so}(\V)$.
Here $\V$ has no cyclic vectors and the previous proofs do not apply. However,
for all possible shapes, there is an argument reducing the
problem to a smaller {\sl rectangle}. \par
1. Consider a near-rectangular diagram of
the first shape:
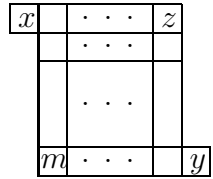
\begin{center}
\begin{picture}(80,70)(0,5) \setlength{\unitlength}{0.015in}
\put(10,10){\line(1,0){60}}
\put(10,20){\line(1,0){60}}
\put(10,50){\line(1,0){50}}
\put(0,60){\line(1,0){60}}
\put(0,70){\line(1,0){60}}
\put(10,10){\line(0,1){60}}
\put(20,10){\line(0,1){60}}
\put(50,10){\line(0,1){60}}
\put(60,10){\line(0,1){60}}
\put(70,10){\line(0,1){10}}
\put(0,60){\line(0,1){10}}
\multiput(25,15)(0,20){3}{. . .}
\put(25,65){. . .}
\put(2,62){$x$}
\put(52,62){$z$}
\put(10,13){$m$}
\put(62,13){$y$}
\end{picture}
\end{center}
Four letters inside of squares stands for the corresponding basis vectors in
$\V$. By assumption, the length of the last and first row is odd, say $2k{+}1$.
Without loss of generality, we may assume that $y=e_1^{2k}m$, $z=e_1^{2k}x$,
and $(x,y)=1$. Then $(z,m)=1$.
\\
Set $\V'=\langle y\rangle^\perp/\langle y\rangle$. Obviously, $\V'$ is again
a quadratic space and $\dim\V'=\dim\V-2$. Put $\g={\frak so}(\V)$,
$\g'={\frak so}(\V')$. Consider
$\p:=\{A\in {\frak so}(\V)\mid Ay\in \langle y\rangle\}$ and the natural map
$\phi:\p\rar {\frak so}(\V')$. Then $\Ker\phi=
\{A\mid A{\cdot}\langle y\rangle^\perp
\subset\langle y\rangle\}$. We are interested in $\Ker\phi\cap\z_\g(\eb)$.
Since $m$ is a cyclic vector in $\langle y\rangle^\perp$ relative to $e_1$
and $e_2$,
we see that any $B\in\Ker\phi\cap\z_\g(\eb)$ takes $m$ to $\langle y\rangle$
and all other basis vectors of $\langle y\rangle^\perp$ to zero.
If $Bm=0$, then $B{\cdot}\langle y\rangle^\perp=0$ and, as $B$ is skew-symmetric,
$B=0$.
Consequently, $\dim(\Ker\phi\cap\z_\g(\eb))\le 1$ \ (actually, $=1$).
\\[1ex]
{\bf Claim:} \ $\z_\g(\eb)\subset\p$. \\ 
Indeed, let $A\in\z_\g(\eb)$. Since $e_iy=0$ for $i=1,2$, we have
$Ay\in\V^{\langle e_1,e_2\rangle}=\langle y,z\rangle$. Now \\
$(Ay,m)=-(y,Am)=-(e_1^{2k}m,Am)=-(m,e_1^{2k}Am)=-(m,Ae_1^{2k}m)=-(m,Ay)$. \\
Whence $(Ay,m)=0$ and therefore $Ay\in\langle y\rangle$.
\\[1ex]
It follows from the claim that $\dim\z_\g(\eb)\le \dim\phi(\z_\g(\eb))+
\dim(\z_\g(\eb)\cap\Ker\phi)=\dim\phi(\z_\g(\eb))+1$. On the other hand,
$\eb$ induces the nilpotent pair $\eb'$ in $\g'$, which corresponds to
the rectangular diagram $\Gamma'=\Gamma\setminus\{x,y\}$. Hence
$\eb'$ is principal in $\g'$. Since $\phi(\z_\g(\eb))\subset\z_{\g'}(\eb')$,
we obtain $\dim\z_\g(\eb)\le\rk\g'+1=\rk\g$. Thus, $\eb$ is principal, too.
\par
2. Up to the permutation of $e_1$ and $e_2$, a near-rectangular diagram of the
second shape is the same a one of the first shape.
\par
3. If $\Gamma$ is a near-rectangular diagram of the third shape, the
same idea can be used. In place of $\langle y\rangle$, one has to consider the
linear span of all basis vectors corresponding to the rightmost column
of $\Gamma$. The details are left to the reader. \par
(ii) We use the notation from the proof of Prop.~\ref{existDd}(ii). \\
(a): If $\Gamma_1$ is rectangular, then $\dim\V_1^{\langle
e_1,e_2\rangle}=1$ (the space $\V_1^{\langle e_1,e_2\rangle}$
corresponds to the sinks of $\Gamma_1$). Next, $\eb$ is principal in
${\frak so}(\V_1)$ by part (i). Whence
$\dim\z_{{\frak so}(\V)}(\eb)=
\rk{\frak so}(\V_1)+1=\rk{\frak so}(\V)$.
\\
(b): Since $\Gamma_1$ has no vertical arrows and
$\Gamma_2$ has no horizontal arrows, we have
$e_2^{(1)}=h_2^{(1)}=0$ and $e_1^{(2)}=h_1^{(2)}=0$. That is,
$e_i$ and $h_i$ operate only in $\V_i$. Moreover, these are regular elements in
${\frak so}(\V_i)$. Now
${\frak so}(\V)\simeq{\frak so}(\V_1)\oplus{\frak so}(\V_2)\oplus
(\V_1\otimes\V_2)$ and $(\V_1\otimes\V_2)^{\langle e_1,e_2\rangle}=
\V_1^{\langle e_1\rangle}\otimes\V_2^{\langle e_2\rangle}$ is
1-dimensional. Therefore
$\dim\z_{{\frak so}(\V)}(\eb)=\rk{\frak so}(\V_1)+
\rk{\frak so}(\V_2)+1=\rk{\frak so}(\V)$. Thus $\eb$ is
a principal pair.
\end{proof}

\sekt{Classification of the principal and distinguished pairs in classical
Lie algebras\nopagebreak}\nopagebreak%
In this section, we prove that the examples of distinguished and principal
pairs in classical Lie algebras described above exhaust all such pairs.
In other words, we attach an appropriate skew-graph to any
distinguished or principal pair in a classical Lie algebra.
The basic reason why this can be achieved
is the existence of simple relationship between the adjoint
and the tautological representation.
\begin{subs} {General scheme of reasoning} \label{general}
\end{subs} Here we give an outline of our approach to
classification. The missing details are supplied in the subsections
devoted to corresponding classical series.
\\[1ex]
Let $\g=\g(\V)$ be a classical Lie algebra, $\V$ being the space of the
tautological representation. Let $\eb$ be a distinguished nilpotent pair
in $\g$ and $\hb$ a characteristics of $\eb$. Our aim is to associate to
$(\eb,\hb)$ a skew-graph. We know
that $\z_\g(\hb)=:\te$ is
a Cartan subalgebra. This is tantamount to saying that ${\Bbb Q}h_1{+}
{\Bbb Q}h_2$ contains a (rational) regular semisimple element, say $h_0$.
Let $\Delta$ be the root system of $\g$ relative to $\te$. Then \\[.5ex]
\hbox to \textwidth{\quad $(\ast)$ \hfil $\ap(h_0)
\ne 0$ for all $\ap\in\Delta$.\hfil}
The key observation is that $(\ast)$ implies that the eigenspaces of
$h_0$ in $\V$ are 1-dimensional, the only exception being the zero weight
for $\g={\frak so}_{2n}$, which can be of multiplicity two.
This will be explained below.
\par
Next, consider the bi-grading $\V=\oplus_{p,q\in {\Bbb Q}}V_{p,q}$
determined by $\hb$.
Since it is a refinement of the grading determined by $h_0$, we have
$\dim\V_{p,q}\le 1$ (exception: $\dim\V_{0,0}\le 2$ for ${\frak so}_{2n}$).
If all the eigenspaces of $\hb$ in $\V$ are 1-dimensional,
we can associate to $\eb$ a skew-graph just by
considering the action of $e_1,e_2$ on these eigenspaces.
Set $n(\Gamma):=\{(p,q)\mid\V_{p,q}\ne 0\}$. By definition,
the arrow $\{(p,q),(p{+}1,q)\}$ is included in $a(\Gamma)$ if and only if
$e_1(\V_{p,q})\ne 0$, and likewise for $e_2$. Then
conditions~\ref{skew}(i),(iii) are obvious, and \ref{skew}(ii) is satisfied,
because $\mbox{tr\,}(h_j)=0$, $j=1,2$. Finally, condition \ref{skew}(iv)
easily follows from the relation $[e_1,e_2]=0$.
Thus, we have attached a skew-graph to any distinguished pair
in ${\frak sl}_{n}$, ${\frak so}_{2n{+}1}$, and ${\frak sp}_{2n}$.
Although $\Gamma$ can be disconnected here, different connected components
have no common nodes.
The structure of connected components will be considered below, for each
series separately.
In the orthogonal and symplectic cases, it follows from the presence of
invariant bilinear form that $\Gamma$ is c.-s.
For ${\frak so}_{2n}$, $\dim\V_{0,0}=0$ or 2. The case $\dim\V_{0,0}=0$ is
covered by the previous argument.
If $\dim\V_{0,0}=2$ (and $\dim\V_{p,q}\le 1$ for $(p,q)\ne (0,0)$), then
a more careful consideration shows that $\V_{0,0}$ splits well with
respect to $e_1$ and $e_2$, and one obtains
a skew-graph with at most three connected c.-s. components,
as in Prop.~\ref{existDd}(iii).
\\[.5ex]
Once we have attached to $(\eb,\hb)$ one of the skew-graphs described
in \re{existA},\re{existB},\re{existC}, and \re{existDd}, it is easily seen
that we can choose a basis of $\V$ so that $e_1,e_2$ will operate on it
according to the respective formulas of sect.\,2.
\\[.5ex]
After describing skew-graphs corresponding to the distinguished pairs,
we turn to classification of principal
pairs. To solve this problem, we exploit the same trick
for all classical series. We show that if $\Gamma$ does not belong to the
appropriate class of skew-graphs, then $\z_\g(\eb)$ contains a `non-positive'
element $x$, i.e. $x\in\z_\g(\eb)_{p,q}$ with $p<0$ or $q<0$. Hence such
$\eb$ can not be principal, see \ref{princ}(iii).
\begin{subs} {Classification for $\g={\frak sl}(\V)$, $\dim\V=n$}
\label{sl}
\end{subs} \\
This classification was already obtained by V.\,Ginzburg. However, our
argument is shorter, because we do not need a description of the whole
space $\z_\g(\eb)$. \\
Let $\esi_1,\dots,\esi_n\in \mbox{Hom\,}(\te,\Bbbk)$ be the
standard weights of $\V$. Then $\Delta=\{\esi_i-\esi_j\mid i\ne j\}$. Hence
$h_0\in\te$ is regular if and only if $\esi_i(h_0)\ne\esi_j(h_0)$ for
$i\ne j$. Hence the eigenspaces of $h_0$ in $\V$ are 1-dimensional.
In \re{general}, we associated to $(\eb,\hb)$ a skew-graph $\Gamma$. Let us
prove that $\Gamma$ is connected.
Assume not and $\Gamma=\Gamma_1\sqcup\Gamma_2$. Let
$\V=\V_1\oplus\V_2$ be the corresponding $(\eb,\hb)$-stable decomposition.
Define $x\in {\frak sl}(\V)$ by $x\vert_{\V_1}=(\dim\V_2)id$,
$x\vert_{\V_2}=-(\dim\V_1)id$. It follows that $x\in\z_\g(\eb)$ is semisimple,
which contradicts the assumption that $\eb$ is distinguished.
Thus, the distinguished pairs in ${\frak sl}(\V)$ are described by
{\sl connected\/} skew-graphs with $\dim\V$ nodes.  \\
 Let us produce a basis of
$\V$ so that $e_1,e_2$ will operate on it according to the formulas in
\re{existA}. This can be done inductively, starting from an arbitrary node
$(i,j)$. Choose $0\ne v_{i,j}\in \V_{i,j}$. Set
$v_{i{+}1,j}:=e_1(v_{i,j})$, $v_{i{-}1,j}:=e_1^{-1}(v_{i,j})\in\V_{i{-}1,j}$
and likewise for $e_2$. Iterating this construction does not lead to a
contradiction, since $[e_1,e_2]=0$.
\\ \indent
It remains to classify the principal pairs. Consider $\Gamma$ as skew-diagram.
We are only concerned with the extreme columns of $\Gamma$.
Let $(i_0,j_0)$ and $(i_0,j_0{+}k)$ (resp. $(i_1,j_1)$ and $(i_1,j_1{+}l)$)
be the coordinates of the extreme squares in the leftmost (resp. rightmost)
column of $\Gamma$, where $k,l\in\Bbb P$. Since $\Gamma$ is a skew-diagram,
$j_0{+}k\ge j_1{+}l$
and $j_0\ge j_1$. Moreover, $\Gamma$ is a {\bf ne}-Young diagram if and only if
$j_0{+}k= j_1{+}l$ and is a {\bf sw}-Young diagram if and only if
$j_0=j_1$. Assume now that $\Gamma$ is not a Young diagram, i.e.,
$j_0{+}k> j_1{+}l$ and $j_0> j_1$.
\par
a) Suppose $k\le l$. Define $x\in{\frak sl}(\V)$ by
$x(v_{i_0,j})=v_{i_1,j_1+l-j_0-k+j}$ for $j_0\le j\le j_0+k$ and
$x(v_{m,n})=0$ for all other basis vectors.
Then $xe_1=e_1x=0$ and therefore $[x,e_1]=0$. It is also easily seen that
$[x,e_2]=0$. Hence $x\in\z_\g(\eb)$. The corresponding eigenvalue is determined
by the shifts of indices in the definition of $x$. The first shift is
$p=i_1-i_0> 0$ and the second shift is $q=j_1+l-j_0-k<0$. Thus, $x$ is a
non-positive element in $\z_\g(\eb)$ and therefore $\eb$ can not
be principal. \par
b)  Suppose $k> l$. Define $x\in{\frak sl}(\V)$ by
$x(v_{i_0,j})=v_{i_1,j_1-j_0+j}$ for $j_0\le j\le j_0+l$ and
$x(v_{m,n})=0$ for all other basis vectors. It is again a non-positive
element in $\z_\g(\eb)$.
\\[.5ex]
We have thus proved that $\eb$ is not principal whenever $\Gamma$ is not a
Young diagram.
\begin{subs} {Classification for $\g={\frak so}(\V)$, $\dim\V=2n{+}1$}
\label{so-odd}
\end{subs}  \\
The weights of $\V$ are $\pm\esi_1,\dots,\pm\esi_n,0$ and the roots are
$\pm\esi_i\pm\esi_j$ $(i\ne j)$, $\pm\esi_i$. If $h_0\in
{\Bbb Q}h_1{+}{\Bbb Q}h_2\subset\te$ is regular,
then $\esi_i(h_0)\ne 0$ and $\esi_i(h_0)\ne\pm\esi_j(h_0)$ $(i\ne j)$.
Hence the eigenspaces of $h_0$ in $\V$ are 1-dimensional.
In \re{general}, we associated to $(\eb,\hb)$ a skew-graph $\Gamma$ with
$2n{+}1$ nodes.
Since 0 is always an eigenvalue of a semisimple element in $\g$, we have
$(0,0)\in n(\Gamma)$. It follows from the presence of
invariant bilinear form on $\V$ that $\Gamma$ is c.-s. Hence the connected
component $\Gamma_0$, containing $(0,0)$, is c.-s., too. Obviously,
$\Gamma_0$ is integral. Suppose $\Gamma$ is disconnected,
$\Gamma=\Gamma_0\sqcup\Gamma_1$. Then $\Gamma_1$ is also c.-s.,
$\# n(\Gamma_1)$ is even, and $(0,0)\not\in n(\Gamma_1)$. Let
$\V=\V_0\oplus\V_1$ be the corresponding $(\eb,\hb)$-stable decomposition.
Clearly, $\eb$ induces a distinguished pair in ${\frak so}(\V_1)$. In the
subsection devoted to the even-dimensional orthogonal case (\ref{so-even}.1),
we show that $\Gamma_1$ has to be connected and non-integral.
Thus, to any distinguished pair in ${\frak so}(\V)$, we have attached a
c.-s. skew-graph having at most two connected components: one is integral and
another (optional) is non-integral.
\\[.5ex]
Construct a basis of $\V$ so that $e_1,e_2$ act on it according to
Eq.\,\re{so-basis}. This is being implemented for each connected component
of $\Gamma$ separately. For $\Gamma_0$, start with $v_{0,0}\in\V_{0,0}$
such that $(v_{0,0},v_{0,0})=1$ and continue inductively. For
$\Gamma_1$, start with $v_{-1/2,-1/2}\in\V_{-1/2,-1/2}$
such that $\Bigl(-e_1(v_{-1/2,-1/2}),e_2(v_{-1/2,-1/2})\Bigr)=1$ and continue
inductively.
\\[.5ex]
It remains to classify principal pairs. If $\Gamma_1\ne\varnothing$, then the
formula $\g\simeq\wedge^2\V$ shows that some eigenvalues
of $\hb$ in $\g$ are fractional. It then follows from \ref{princ}(ii) that
$\Gamma$ has to be connected whenever $\eb$ is principal. Until the end of
this subsection we assume that $\Gamma$ is connected and integral.
Let $(-i_0,-j_0)$ and $(-i_0,-j_0{+}k)$ (resp. $(i_0,j_0-k)$ and $(i_0,j_0)$)
be the coordinates of the extreme squares in the leftmost (resp. rightmost)
column of $\Gamma$, see Fig.~\ref{extreme}.
\begin{figure}[ht]
\caption{}  \label{extreme}
\setlength{\unitlength}{0.015in}
\centerline
{
\begin{picture}(120,100)
\put(60,0){\line(1,0){50}}
\put(30,10){\line(1,0){30}}
\put(100,10){\line(1,0){10}}
\put(0,20){\line(1,0){30}}
\put(0,30){\line(1,0){10}}
\put(100,50){\line(1,0){10}}
\put(80,60){\line(1,0){30}}
\put(0,70){\line(1,0){10}}
\put(50,70){\line(1,0){30}}
\put(0,80){\line(1,0){50}}
\put(0,20){\line(0,1){60}}
\put(10,20){\line(0,1){60}}
\put(30,10){\line(0,1){10}}
\put(50,70){\line(0,1){10}}
\put(60,0){\line(0,1){10}}
\put(80,60){\line(0,1){10}}
\put(100,0){\line(0,1){60}}
\put(110,0){\line(0,1){60}}
\put(-10,25){\vector(1,0){15}}
\put(-55,24){{\footnotesize $(-i_0,-j_0)$}}
\put(-10,75){\vector(1,0){15}}
\put(-65,74){{\footnotesize $(-i_0,-j_0{+}k)$}}
\put(120,5){\vector(-1,0){15}}
\put(123,4){{\footnotesize $(i_0,j_0{-}k)$}}
\put(120,55){\vector(-1,0){15}}
\put(123,54){{\footnotesize $(i_0,j_0)$}}
\end{picture}
}
\end{figure}

\noindent Since $\Gamma$ is a skew-diagram, we
have $-j_0{+}k\ge j_0$. Being c.-s., $\Gamma$ is a rectangle if and only if
$-j_0{+}k=j_0$, i.e.
$k= 2j_0$. Maintain the convention of \re{existB} concerning the symmetric
form, $e_1$, and $e_2$.
Assuming that $k>2j_0$, we can construct a non-positive element in
$\z_\g(\eb)$ as follows.
\par
a) Suppose $k$ is odd. Define $x\in \mbox{End\,}(\V)$ by
$x(v_{-i_0,j})=(-1)^{j}v_{i_0,2j_0-k{+}j}$ ($-j_0\le j\le -j_0{+}k$) and
$x(v_{m,n})=0$ for all other basis vectors.
Using the formulae from \re{existB}, it is not hard to verify that
$x\in {\frak so}(\V)$ and $[x,e_1]=[x,e_2]=0$.
Thus $x\in\z_\g(\eb)_{p,q}$, where $p=2i_0$ and $q=2j_0-k <0$.
\par
b) Suppose $k$ is even and hence $k>2j_0{+}1$. If $k\ge 2$,
define $x\in {\frak so}(\V)$ by
$x(v_{-i_0,j})=(-1)^{j}v_{i_0,2j_0{+}1-k{+}j}$ ($-j_0\le j\le -j_0{+}k-1$) and
$x(v_{m,n})=0$ for all other basis vectors. Then $x\in\z_\g(\eb)_{p,q}$,
where $p=2i_0$ and $q=2j_0{+1}-k <0$. \\
If $k=0$, then the extreme columns consist of a single square. In this case,
the previous argument can be applied to the extreme {\sl rows\/} of $\Gamma$.
Indeed, as $\Gamma$ is connected, the extreme rows can not consist of a
single square.
\\[.5ex]
We have thus proved that $\eb$ is not principal whenever $\Gamma$ is not a
rectangle.
\begin{subs} {Classification for $\g={\frak sp}(\V)$, $\dim\V=2n$}
\label{sp}
\end{subs} \\
The weights of $\V$ are $\pm\esi_1,\dots,\pm\esi_n$ and the roots are
$\pm\esi_i\pm\esi_j$ $(i\ne j)$, $\pm 2\esi_i$. If $h_0\in
{\Bbb Q}h_1{+}{\Bbb Q}h_2\subset\te$ is regular,
then $\esi_i(h_0)\ne 0$ and $\esi_i(h_0)\ne\pm\esi_j(h_0)$ $(i\ne j)$.
Hence the eigenvalues of $h_0$ in $\V$ are simple and nonzero.
In \re{general}, we  associated to $(\eb,\hb)$ a skew-graph $\Gamma$ with
$2n$ nodes.
Since $0$ is not an eigenvalue of $h_0$, we have $(0,0)\not\in n(\Gamma)$.
For the same reason as in \re{so-odd}, $\Gamma$ is c.-s.
Let $\Gamma_1$ be a connected component of $\Gamma$. Assume that $\Gamma_1$
is not c.-s. Then $-n(\Gamma_1)$ is the set of nodes of another connected
component, say $\Gamma_2$. Therefore there is the decomposition
$\Gamma=\Gamma_1\sqcup\Gamma_2
\sqcup\Gamma_3$, where $\Gamma_3$ is c.-s.
(possibly, $\Gamma_3=\varnothing$).
If $\V=\V_1\oplus\V_2\oplus\V_3$ is the corresponding $(\eb,\hb)$-stable
decomposition, then define $x\in {\frak sp}(\V)$  by
$x\vert_{\V_1}=id$, $x\vert_{\V_2}=-id$, $x\vert_{\V_3}=0$.
Then $x$ is a semisimple element in $\z_\g(\eb)$, which contradicts the fact
that $\eb$ is distinguished. Hence each connected component must be
c.-s. Let $\Gamma_\ap$ be one of them. Since $(0,0)\not\in n(\Gamma_\ap)$,
$\Gamma_\ap$ is not integral and we may assume without loss that $j\in
\frac{1}{2}+\Bbb Z$ for any $(i,j)\in n(\Gamma_\ap)$.
\\[.5ex]
{\bf Lemma.} {\sl Then} $i\in\Bbb Z$. \\
{\sl Proof.} Let $(-m,-n)$ be a node of $\Gamma_\ap$ in the negative quadrant
and $v_{-m,-n}\in\V_{-m,-n}$ a corresponding weight vector.
Such a node exists because $\Gamma_\ap$ is c.-s. and connected. Then
$0\ne e_1^{2m}e_2^{2n}v_{-m,-n}\in \V_{m,n}$.
By assumption, $2n$ is odd. Assume that $2m$ is odd, too. Using the
{\sl alternating\/} bilinear form on $\V$, we then obtain
\[
(e_1^{2m}e_2^{2n}v_{-m,-n},v_{-m,-n})=(v_{-m,-n},
e_1^{2m}e_2^{2n}v_{-m,-n})=-(e_1^{2m}e_2^{2n}v_{-m,-n},v_{-m,-n}) \ .
\]
Whence $(e_1^{2m}e_2^{2n}v_{-m,-n},v_{-m,-n})=0$, which
contradicts the fact that $(~,~)$ is non-degenerate. Thus,
$2m$ must be even. \qus
The lemma means $\Gamma_\ap$ is semi-integral. Recall that there are two
sorts of semi-integral c.-s. diagrams. Two diagrams of different sorts
have disjoint sets of nodes, while diagrams of the same sort always have
some common nodes. Since
$n(\Gamma_\ap)\cap n(\Gamma_\beta)=\varnothing$ for any couple of connected
components, we obtain $\Gamma$ has at most two connected components,
necessarily of different sorts.
\\[.6ex]
Constructing a basis of $\V$ so that $e_1,e_2$ act on it according to
Eq.\,\re{sp-basis} starts with choosing $v_{0,-1/2}\in\V_{0,-1/2}$
so that $\Bigl(e_2(v_{0,-1/2}),v_{0,-1/2}\Bigr)=1$. For the connected
component of another sort, start with $v_{-1/2,0}$ and so on..
\\[.6ex]
It remains to describe the skew-graphs arising from principal pairs.
If $\Gamma$ has two connected components (of different sorts), then the
formula $\g\simeq {\cal S}^2\V$ shows that some
eigenvalues of $\hb$ in $\g$ are fractional. It then follows from
\ref{princ}(ii) that
$\Gamma$ has to be connected whenever $\eb$ is principal.
Without loss, we may assume that $i\in\Bbb Z$ for any $(i,j)\in n(\Gamma)$.
Again, we only consider the extreme columns of $\Gamma$,
see Figure~\ref{extreme}. (If $\Gamma$ is
semi-integral of another sort, one has to consider extreme rows and perform
`transposed' constructions.)
As above,  $-j_0{+}k\ge j_0$, and $\Gamma$ is a rectangle if and only if
$k= 2j_0$, which is now odd\,!
Maintain the conventions of \re{existC} concerning the alternating form,
$e_1$, and $e_2$.
Assuming that $k>2j_0$, we can construct a non-positive element in
$\z_\g(\eb)$ as follows:
\par
a) Suppose $k$ is even. Define the endomorphism  $x$ by 
$x(v_{-i_0,j})=\ap v_{i_0,2j_0{-}k{+}j}$ and $x(v_{m,n})=0$ for all
other basis vectors, where
$\ap=\left\{\begin{array}{ll} 1, & \mbox{if } \ -j_0\le j<0 \\
(-1)^{i_0{+}j{+}\frac{1}{2}}, & \mbox{if } \ 0< j< k{-}2j_0 \\
-1, & \mbox{if } \ k{-}2j_0< j\le k{-}j_0
\end{array}.\right.$
A straightforward verification shows that $x$ is symplectic and
$x\in\z_\g(\eb)_{p,q}$ with
$q=2j_0-k<0$.
\par
b) Suppose $k$ is odd and hence $k>2j_0{+}1$. Define $x$ by
$x(v_{-i_0,j})=\ap v_{i_0,2j_0{-}k{+}j{+}1}$ and $x(v_{m,n})=0$ for all
other basis vectors, where \\
$\ap=\left\{\begin{array}{ll} 1, & \mbox{if} \ -j_0\le j<0 \\
(-1)^{i_0{+}j{+}\frac{1}{2}}, & \mbox{if} \ 0< j< k{-}2j_0{-}1 \\
-1, & \mbox{if} \ k{-}2j_0{-}1< j\le k{-}j_0{-}1
\end{array} .\right.$
A straightforward verification shows that  $x$ is symplectic and
$x\in\z_\g(\eb)_{p,q}$ with
$q=2j_0{+}1-k<0$.
\\[.6ex]
We have thus proved that $\eb$ is not principal whenever $\Gamma$ is not
a rectangle.
\begin{subs} {Classification for $\g={\frak so}(\V)$, $\dim\V=2n$}
\label{so-even}
\end{subs} \\
The weights of $\V$ are $\pm\esi_1,\dots,\pm\esi_n$ and the roots are
$\pm\esi_i\pm\esi_j$ $(i\ne j)$. If $h_0\in
{\Bbb Q}h_1{+}{\Bbb Q}h_2\subset\te$ is regular,
then $\esi_i(h_0)\ne\pm\esi_j(h_0)$ $(i\ne j)$.
Hence there are two possibilities: either all $\esi_i(h_0)$ are non-zero,
or $\esi_{l}(h_0)=0$ for a unique $l\in\{1,\dots,n\}$.
Consider these cases in turn.
\\[1ex]
{\bf (\theequation.1)} \quad  $\esi_i(h_0)\ne 0$ for all $i$. \\
Here $\V_{0,0}=0$ and $\dim\V_{p,q}\le 1$ for all $(p,q)$. Hence the argument
from \re{general} associates to such $\eb$ a skew-graph $\Gamma$ with
$2n$ nodes. It follows from the presence of invariant bilinear form on $\V$
that $\Gamma$ is c.-s. Since $0$ is not an eigenvalue of $h_0$, we have
$(0,0)\not\in n(\Gamma)$. The same argument as in \re{sp} proves that each
connected component of $\Gamma$ is also c.s. (otherwise we would find
a semisimple $x\in\z_\g(\eb)$).
Let $\Gamma_\ap$ be one of them. Since $(0,0)\not\in n(\Gamma_\ap)$,
$\Gamma_\ap$ is not integral and we may assume without loss that $j\in
\frac{1}{2}+\Bbb Z$ for any $(i,j)\in n(\Gamma_\ap)$.
\\[.5ex]
{\bf Lemma.} {\sl Then} $i\in\frac{1}{2}+\Bbb Z$, too. \\
{\sl Proof.} Let $(-m,-n)$ be a node of $\Gamma_\ap$ in the negative quadrant
and $v_{-m,-n}\in\V_{-m,-n}$ a corresponding weight vector.
Such a node exists because $\Gamma_\ap$ is c.-s. and connected. Then
$0\ne e_1^{2m}e_2^{2n}v_{-m,-n}\in \V_{m,n}$.
By assumption, $2n$ is odd.
Assume that $2m$ is even. Using the {\sl symmetric\/} bilinear form on $\V$, we
then obtain
\[
(e_1^{2m}e_2^{2n}v_{-m,-n},v_{-m,-n})=-(v_{-m,-n},
e_1^{2m}e_2^{2n}v_{-m,-n})=-(e_1^{2m}e_2^{2n}v_{-m,-n},v_{-m,-n}) \ .
\]
Whence $(e_1^{2m}e_2^{2n}v_{-m,-n},v_{-m,-n})=0$, which
contradicts the fact that $(~,~)$ is non-degenerate. Thus,
$2m$ must be odd. \qus
The lemma means that each $\Gamma_\ap$ is non-integral. Obviously, any
connected non-integral c.-s. skew-graph contains the node $(1/2,1/2)$.
Whence $\Gamma$ must be connected.
Let us realize which skew-graphs arise in connection with principal pairs.
\\[1ex]
I. As in \re{so-odd} and \re{sp}, consider the extreme columns of $\Gamma$
and their extreme squares (Figure~\ref{extreme}). In this case,
$i_0,j_0\in\frac{1}{2}+\Bbb Z$. The non-negative integer $k-2j_0$ represents
the relative {\sl vertical\/} shift of the two columns.
\par
a) \un{$k$ is odd}. Then $k-2j_0$ is even. Suppose $k-2j_0\ge 2$. Define the
endomorphism $x$ by
$x(v_{-i_0,j})=(-1)^{j}v_{i_0,2j_0-k{+}j}$ ($-j_0\le j\le -j_0{+}k$) and
$x(v_{m,n})=0$ for all other basis vectors. Then $x\in{\frak so}(\V)$ and
$x\in\z_\g(\eb)_{p,q}$ with $q=2j_0-k <0$. \par
Thus, if $\eb$ is principal and $k$ is odd, then the vertical shift is
0, i.e., $\Gamma$ is rectangular.
\par
a) \un{$k$ is even}. Then $k-2j_0$ is odd. Suppose $k-2j_0\ge 3$ and $k>0$.
Then $x$ can be defined  by
$x(v_{-i_0,j})=(-1)^{j}v_{i_0,2j_0-k{+}j{+}1}$ ($-j_0\le j\le -j_0{+}k-1$) and
$x(v_{m,n})=0$ for all other basis vectors. Here
$x\in\z_\g(\eb)_{p,q}$ with $q=2j_0{-}k{+}1 <0$.  \par
Thus, if  $\eb$ is principal and $k$ is even, then either the vertical shift is
equal to 1 or $k=0$, i.e., the extreme columns consist of a sole square.
\\[1ex]
II. The same argument applies to the extreme rows of $\Gamma$, and
we obtain the same conditions for their length and possible
{\sl horizontal\/} shift.
\\[1ex]
III. Combining these constraints together yield exactly the class of
near-rectangular diagrams, see \re{near}. It has to only be noted that one
combination can not occur.
Namely, if the length of all extreme rows and columns were
equal to 1, then $\Gamma$ appeared to be disconnected.
\\[1ex]
{\bf (\theequation.2)} \quad There is a unique $l$ such that $\esi_l(h_0)=0$.\\
Then $\dim\V_{0,0}=0$ or 2, because $\dim\V$ is even and $\dim\V_{p,q}=
\dim\V_{-p,-q}$. The case $\dim\V_{0,0}=0$ belongs to the previous
part. Let $\dim\V_{0,0}=2$.
Our aim is to associate to such $(\eb,\hb)$ a skew-graph with at most
three connected components, as in Prop.~\ref{existDd}(iii).
Let $\V_{int}$ (resp. $\V_{f\! r}$) be the sum of all integral (resp.
fractional) eigenspaces of $\hb$ in $\V$ (see sect.~1 for this terminology).
Clearly, $\V_{int}\perp \V_{f\! r}$. Then the argument from (\ref{so-even}.1)
applies to $\V_{f\! r}$ and yields a connected c.-s. non-integral skew-graph
$\Gamma_0$. If $\V_{f\! r}=0$, then $\Gamma_0=\varnothing$. However,
$\V_{int}\ne 0$. Although constructing a skew-graph associated to
$\V_{int}$ contains no extraordinary ideas,
it is rather long and tedious. For this reason, we only sketch the main steps.
The following argument takes place inside of $\V_{int}$, i.e.,
all indices are assumed to be integral.\\
{\bf I.} Assume that $\V_{-p,-q}\ne 0$ for some $p,q>0$. \\
Consider the minimal
subspace $W\supset \V_{-p,-q}$ satisfying the properties: \par
1. $W$ is $(\eb,\hb)$-stable; \par
2. if $x\in\V$ is an $\hb$-eigenvector, $x\not\in \V_{0,0}$, and
$e_i(x)\subset W$ for $i=1$ or 2, then $x\in W$. \\
Let us say that $W$ is the {\it envelope\/} of $\V_{-p,-q}$.
It is easily seen that $
\bigoplus_{i\ge -p,\,j\ge -q}
W_{i,j}$ is generated by $W_{-p,-q}=\V_{-p,-q}$ as $\eb$-module. In particular,
$\dim W_{0,0}\le 1$.
This already means that one can associate to $W$ a (connected\,!) graph.
Consider an $\hb$-stable decomposition
$W=\tilde W\oplus (W\cap W^\perp)$.
Using the definition of $W$, the fact that $W\cap W^\perp$ is also
$\eb$-stable, and that $(~,~)\vert_{\tilde W}$ is non-degenerate,
one proves that both $\tilde W$ and $W\cap W^\perp$ satisfy
properties 1 and 2 above. Whence either $\tilde W=0$ or $W\cap W^\perp=0$.
If the first relation were true, we would find a semisimple element
$x\in\z_\g(\eb)$ (cf. \ref{sp}). Thus, $W\cap W^\perp=0$ and
$(~,~)\vert_{W}$ is non-degenerate. Consequently, $W_{\ap,\beta}\ne 0$
if and only if $W_{-\ap,-\beta}\ne 0$. Since $W_{-p,-q}\ne 0$, it follows
from the definition of $W$ that $W_{p,q}=e_1^{2p}e_2^{2q}(W_{-p,-q})\ne 0$.
In particular, $\V_{0,0}=W_{0,0}=e_1^{p}e_2^{q}(W_{-p,-q})\ne 0$.
Now, it is easily seen that $W$ does not
depend on initial choice of $(-p,-q)$. Therefore $\V_{\ap,\beta}\subset W$
for all $(\ap,\beta)$ with $\ap\beta>0$. This implies that $W_{p,q}
=\V_{p,q}$ for all $(p,q)\ne (0,0)$. (Otherwise we would produce a
semisimple $x\in\z_\g(\eb)$, cf. \ref{sp}.)
Thus, we obtain the $(\eb,\hb)$-stable decomposition
$\V_{int}=W\oplus\{\mbox{1-dim space}\}$, where $\V_{0,0}=W_{0,0}
\oplus\{\mbox{1-dim space}\}$. Clearly,
this situation corresponds to a skew-graph described in
Prop.~\ref{existDd}(ii)a. \\
{\bf II.} Assume now that $\V_{p,q}=0$ for all $(p,q)$ with $pq>0$. \\
If $\V_{p,0}\ne 0$ for some $p<0$, then consider the envelope of
$\V_{p,0}$, as in part I. Call it $W_1$. Similarly, if $\V_{0,q}\ne 0$
for some $q<0$, consider $W_2$, the envelope of $\V_{0,q}$.
As in part I, one proves that
$(~,~)\vert_{W_i}$ is non-degenerate ($i=1,2$). Now, one of the two things can
occur: either $W_1=W_2$ and this space is of codimension 1 in $\V$, or
$W_1\oplus W_2=\V$. The first possibility still belongs to the case studied
in part I and corresponds to a skew-graph described in \ref{existDd}(ii)a,
whereas the second one corresponds to  \ref{existDd}(ii)b.
\\[.6ex]
The case, where $\V_{p,0}=0$ and $\V_{0,q}=0$ for all $p,q\ne 0$,
is impossible. For, otherwise the $(\eb,\hb)$-stable
decomposition $\V_{int}=\V'\oplus\V''\oplus\V_{0,0}$ with isotropic spaces
$\V'$ and $\V''$ would produce a semisimple element $x\in\z_\g(\eb)$.
(Set $x\vert_{\V'}=id$, $x\vert_{\V''}=-id$, and $x\vert_{\V_{0,0}}=0$.)
\\[.6ex]
Finally, if $\V_{p,0}\ne 0$ for some $p\ne 0$ but
$\V_{0,q}=0$ for all $q\ne 0$, then the envelope $W_1$ still exists. Here
one proves that $W_1\oplus\{\mbox{1-dim space}\}=\V_{int}$. Again,
this corresponds to \ref{existDd}(ii)a.
\\[.5ex]
Thus, any distinguished pair with $\dim\V_{0,0}=2$ can be described via
a skew-graph with at most three connected components, as in \ref{existDd}(iii).
Because connected components of $\Gamma$ are of the same type as in the
odd-dimensional orthogonal case, the procedure of constructing
basis of $\V$ so that $e_1,e_2$ act on it according to Eq.\,\re{so-basis}
is the same as in  \ref{so-odd}.
\\[.5ex]
It remains to describe the skew-graphs
arising from principal pairs. We have attached to $(\eb,\hb)$ a c.-s.
skew-graph $\Gamma_0\cup\Gamma_1\cup\Gamma_2$, where $\Gamma_0$
is non-integral and $\Gamma_1, \Gamma_2$ are integral. If $\Gamma_0\ne
\varnothing$, then $\hb$ has fractional eigenvalues in $\g=\wedge^2\V$.
Hence $\Gamma_0=\varnothing$ for principal pairs. We have now two
possibilities for $\Gamma_1$ and $\Gamma_2$, described in \ref{existDd}(ii).
\par
Consider case a). Here $\dim\V_2=1$, $\dim\V_1$ is odd, and
$\g={\frak so}(\V)\simeq {\frak so}(\V_1)\oplus \V_1$. Therefore
$\z_{\g}(\eb)={\frak so}(\V_1)^{\langle e_1,e_2\rangle}
\oplus \V_1^{\langle e_1,e_2\rangle}=\z_{{\frak so}(\V_1)}(\eb)\oplus
\V_1^{\langle e_1,e_2\rangle}$. Whence $\dim\z_\g(\eb)\ge \rk{\frak so}(\V_1)
{+}1=\rk {\frak so}(\V)$, and the equality holds if and only if $\eb$
is a principal pair in ${\frak so}(\V_1)$ and $\dim
\V_1^{\langle e_1,e_2\rangle}=1$. Each of the two conditions implies
that $\Gamma_1$ is rectangular (for the first one see \ref{so-odd}).
\par
Consider case b).
Here $\z_{\g}(\eb)={\frak so}(\V_1)^{\langle e_1,e_2\rangle}
\oplus{\frak so}(\V_2)^{\langle e_1,e_2\rangle}
\oplus (\V_1\otimes\V_2)^{\langle e_1,e_2\rangle}$.
Hence $\eb$ is principal in $\g$ if and only if $\eb$ is principal in
${\frak so}(\V_i)$ ($i=1,2$) and
$\dim(\V_1\otimes\V_2)^{\langle e_1,e_2\rangle}=1$. By \ref{so-odd}, the first
two conditions imply that $\Gamma_1$ and $\Gamma_2$ are rectangular.
Moreover, as $n(\Gamma_1)\cap n(\Gamma_2)=\{(0,0)\}$, $\Gamma_1$ has to be
a horizontal chain and $\Gamma_1$ has to be
a vertical chain. This means that $e_i$ operate only in $\V_i$ and $e_i$
is a regular nilpotent element in ${\frak so}(\V_i)$.
Then
$(\V_1\otimes\V_2)^{\langle e_1,e_2\rangle}=
\V_1^{\langle e_1\rangle}\otimes\V_2^{\langle e_2\rangle}$
is 1-dimensional. \\
We have thus proved that the skew-graphs associated to principal pairs
are those described in Prop.~\ref{existDp}(ii).
\\[1ex]
Thus, the classification of the distinguished and principal nilpotent pairs
in the classical simple Lie algebras is completed.
We summarize it in the following
\begin{s}{THEOREM}  \label{MAIN}
Let $\g$ be a classical simple Lie algebra. \par
{\sf 1.} There is a mapping from
the set of $G$-orbits of distinguished pair in $\g$ to the
appropriate set of skew-graphs:
$G{\cdot}\eb\stackrel{\vp}{\mapsto} \Gamma(\eb)$.
This mapping is a bijection, unless  $\g={\frak so}_{2n}$ and
$\Gamma(\eb)$ is connected. In the last case, $\vp^{-1}(\Gamma(\eb))$
consists of two $SO_{2n}$-orbits.
The description of the relevant skew-graphs is given in
\ref{existA}(1),\,\ref{existB}(1),\,\ref{existC}(1), and \ref{existDd},
while constructing $\vp$ is accomplished in \re{general}--\re{so-even}.
\par
{\sf 2.} The statement with the word `principal' in place of `distinguished'
is also true. In this case the description of the relevant skew-graphs is found
in \ref{existA}(2),\,\ref{existB}(2),\,\ref{existC}(2), and \ref{existDp}.
\end{s}%
Since the number of relevant skew-graphs with a fixed cardinality of
$n(\Gamma)$ is finite, the number of $G$-orbits of the distinguished
pairs is finite, too.

\sekt{Some complements to the classification\nopagebreak}%
\begin{subs}{Explicit description of the centralizer}
\end{subs}
Here we give a basis of $\z_\g(\eb)$ for all principal nilpotent pairs.
This work was partly fulfilled in sect.\,2.
Now we complete it and present the answer systematically. In the following
formulas $k,l\in\Bbb P$. Note that $e_1^k e_2^l\in\z_\g(\eb)_{k,l}$ whenever
$e_1^k e_2^l\in\g$.
\begin{subsubs}{$\g={\frak sl}(\V)$} \end{subsubs} \\
Let $\Gamma$ be a {\bf sw}-Young graph and $(i_0,j_0)$ the source of it.
Then
\[
\z_\g(\eb)=\langle e_1^k e_2^l \mid (i_0+k,j_0+l)\in n(\Gamma) \ \&
\ (k,l)\ne (0,0) \rangle \ .
\]
Let $\Gamma$ be a {\bf ne}-Young graph and $(i_1,j_1)$ the sink of it.
Then
\[
\z_\g(\eb)=\langle e_1^k e_2^l \mid (i_1-k,j_1-l)\in n(\Gamma) \ \&
\ (k,l)\ne (0,0) \rangle \ .
\]
\begin{subsubs}{$\g={\frak so}(\V)$ or ${\frak sp}(\V)$}  \end{subsubs} \\
Let $\Gamma$ be a connected rectangular skew-graph (of appropriate type) and
$(i_0,j_0)$ the source of it. Then
\[
\z_\g(\eb)=\langle e_1^k e_2^l \mid (i_0+k,j_0+l)\in n(\Gamma) \ \&
\ k{+}l \mbox{ is odd } \rangle \ .
\]
\begin{subsubs}{$\g={\frak so}(\V)$, $\dim\V$ is even}  \end{subsubs}
\\
In the following figures, letters $x,y,z,u$ inside of some squares denote
corresponding basis vectors in $\V$. It is assumed that $(x,y)=(u,z)=1$.
Describing operator $A\in{\frak so}(\V)$, we indicate only its nonzero entries.
\\[.5ex]
$\bullet$ Let $\Gamma$ be a near-rectangular diagram of the first
shape (see \re{near}):
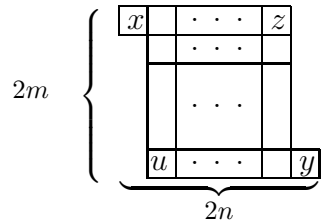
\begin{center}
\setlength{\unitlength}{0.015in}
\raisebox{36\unitlength}%
{$\left\{\parbox{1pt}{\vspace{60\unitlength}}\right.$}
\begin{picture}(80,77)
\put(-38,38){{\footnotesize $2m$}}
\put(10,10){\line(1,0){60}}
\put(10,20){\line(1,0){60}}
\put(10,50){\line(1,0){50}}
\put(0,60){\line(1,0){60}}
\put(0,70){\line(1,0){60}}
\put(10,10){\line(0,1){60}}
\put(20,10){\line(0,1){60}}
\put(50,10){\line(0,1){60}}
\put(60,10){\line(0,1){60}}
\put(70,10){\line(0,1){10}}
\put(0,60){\line(0,1){10}}
\multiput(25,15)(0,20){3}{. . .}
\put(25,65){. . .}
\put(2,62){$x$}
\put(52,62){$z$}
\put(10,13){$u$}
\put(62,13){$y$}
\put(-.5,10){$\underbrace%
{\mbox{\hspace{69\unitlength}}}_{\mbox{{\footnotesize $2n$}}}$}
\end{picture}
\end{center}
Then $\rk\g=2m(n-1)+1$ and
\[
\z_\g(\eb)=\langle e_1^k e_2^l \mid k\le 2n-3,\ l\le 2m-1 \ \& \
\ k{+}l \mbox{ is odd } \rangle\oplus \langle A\rangle \ ,
\]
where $Ax=z$ and $Au=-y$. Here $A\in\z_\g(\eb)_{2n-2,0}$.
\\
$\bullet$ The description for the second shape is the `transposed' one.
\\
$\bullet$ Let $\Gamma$ be a near-rectangular diagram of the third shape
(see \re{near}):
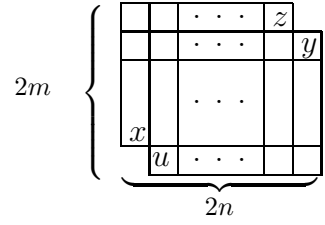
\begin{center}
\setlength{\unitlength}{0.015in}
\raisebox{36\unitlength}%
{$\left\{\parbox{1pt}{\vspace{60\unitlength}}\right.$}
\begin{picture}(80,77)
\put(-38,38){{\footnotesize $2m$}}
\put(10,10){\line(1,0){60}}
\put(0,20){\line(1,0){70}}
\put(0,50){\line(1,0){70}}
\put(0,60){\line(1,0){70}}
\put(0,70){\line(1,0){60}}
\put(10,10){\line(0,1){60}}
\put(20,10){\line(0,1){60}}
\put(50,10){\line(0,1){60}}
\put(60,10){\line(0,1){60}}
\put(70,10){\line(0,1){50}}
\put(0,20){\line(0,1){50}}
\multiput(25,15)(0,20){3}{. . .}
\put(25,65){. . .}
\put(2,22){$x$}
\put(52,62){$z$}
\put(10,13){$u$}
\put(62,53){$y$}
\put(-.5,10){$\underbrace%
{\mbox{\hspace{68\unitlength}}}_{\mbox{{\footnotesize $2n$}}}$}
\end{picture}
\end{center}
Then $\rk\g=2mn-1$ and
\[
\z_\g(\eb)=\langle e_1^k e_2^l \mid k\le 2n-1,\ l\le 2m-1 \ \& \
\ k{+}l \mbox{ is odd }\& \ k+l\ne 2m{+}2n{-}3 \rangle
\oplus \langle A\rangle \ ,
\]
where $Ax=z$ and $Au=-y$. Here $A\in\z_\g(\eb)_{2n-2,2m-2}$.
\\
$\bullet$ Let $\Gamma$ be disconnected, as in \ref{existDp}(ii)b:
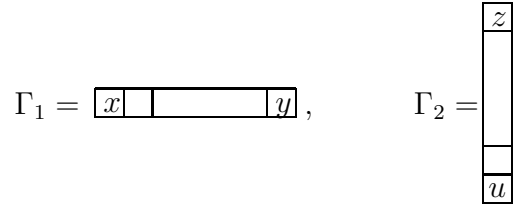
\begin{center}
\setlength{\unitlength}{0.015in}
\begin{picture}(80,70)
\put(-29,31){$\Gamma_1=$}
\put(0,30){\line(1,0){70}}
\put(0,40){\line(1,0){70}}
\put(0,30){\line(0,1){10}}
\put(10,30){\line(0,1){10}}
\put(20,30){\line(0,1){10}}
\put(60,30){\line(0,1){10}}
\put(70,30){\line(0,1){10}}
\put(2,32){$x$}
\put(62,32){$y$}
\put(73,31){,}
\end{picture}
\qquad
\begin{picture}(40,70)
\put(5,31){$\Gamma_2=$}
\put(30,0){\line(0,1){70}}
\put(40,0){\line(0,1){70}}
\put(30,0){\line(1,0){10}}
\put(30,10){\line(1,0){10}}
\put(30,20){\line(1,0){10}}
\put(30,60){\line(1,0){10}}
\put(30,70){\line(1,0){10}}
\put(31,2){$u$}
\put(32,62){$z$}
\end{picture}
\end{center}
Assume that $\# n(\Gamma_1)=2n{+}1$ and $\# n(\Gamma_2)=2m{+}1$. Then
$\rk\g=n{+}m{+}1$ and
\[
\z_\g(\eb)=\langle e_1,e_1^3,\dots,e_1^{2n-1},e_2,e_2^3,\dots,e_2^{2m-1}
\rangle\oplus\langle A\rangle \ ,
\]
where $Ax=z$ and $Au=-y$. Here $A\in\z_\g(\eb)_{n,m}$.
\\
$\bullet$ Let $\Gamma$ be disconnected, as in \ref{existDp}(ii)a:
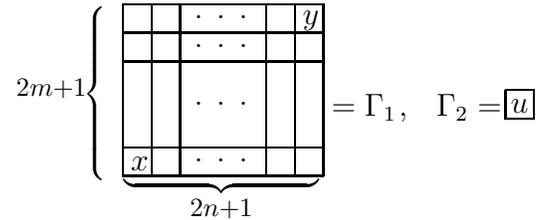
\begin{center}
\setlength{\unitlength}{0.015in}
\raisebox{36\unitlength}%
{$\left\{\parbox{1pt}{\vspace{60\unitlength}}\right.$}
\begin{picture}(110,70)
\put(-38,38){{\footnotesize $2m{+}1$}}
\put(72,31){$=\Gamma_1$,}
\put(0,10){\line(1,0){70}}\put(0,20){\line(1,0){70}}
\put(0,50){\line(1,0){70}}\put(0,60){\line(1,0){70}}
\put(0,70){\line(1,0){70}}
\put(10,10){\line(0,1){60}}\put(20,10){\line(0,1){60}}
\put(50,10){\line(0,1){60}}\put(60,10){\line(0,1){60}}
\put(70,10){\line(0,1){60}}\put(0,10){\line(0,1){60}}
\multiput(25,15)(0,20){3}{. . .}
\put(25,65){. . .}
\put(2,12){$x$}
\put(62,63){$y$}
\put(-.5,10){$\underbrace%
{\mbox{\hspace{68\unitlength}}}_{\mbox{{\footnotesize $2n{+}1$}}}$}
\end{picture}
\begin{picture}(30,70)
\put(-5,31){$\Gamma_2=$}
\put(20,30){\line(1,0){10}}\put(20,40){\line(1,0){10}}
\put(20,30){\line(0,1){10}}\put(30,30){\line(0,1){10}}
\put(21,32){$u$}
\end{picture}
\end{center}
Now $u=z$, i.e., $(u,u)=1$. Then $\rk\g=2mn{+}m{+}n{+}1$ and
\[
\z_\g(\eb)=\langle e_1^k e_2^l \mid k\le 2n,\ l\le 2m \ \& \
\ k{+}l \mbox{ is odd }\rangle
\oplus \langle A\rangle \ ,
\]
where $Ax=u$ and $Au=-y$. Here $A\in\z_\g(\eb)_{n,m}$.
\\[.6ex]
In all cases, it is clear that $A\in\z_\g(\eb)$. Some explanation is only
needed for the eigenvalue of $A$. Recall that all the diagrams above are
endowed with the coordinate system so that the origin is in the barycentre of
each of them. Then the coordinate of the centre of each square is the
eigenvalue of the corresponding $\hb$-eigenspace of $\V$. In the last
case, we have $y\in\V_{n,m}$, $x\in\V_{-n,-m}$, and
$u\in\V_{0,0}$. Therefore $[h_1,A]=nA$ and $[h_2,A]=mA$.
\\[.6ex]
In \cite{vitya}, the eigenvalues of $\hb$ in $\z_\g(\eb)$ were called
the {\it bi-exponents\/} corresponding to $\eb$. These are important for
some application of principal pairs. Thus, we have also described the
bi-exponents for all principal pairs.

\begin{subs}{On rectangular nilpotent pairs}
\end{subs}%
In \cite{wir}, we introduced the notion of a rectangular principal pair and
gave a classification of such pairs. The fact that for ${\frak sl}(\V)$
a principal pair is rectangular if and only if $\Gamma$ is rectangular
was the only motivation for the name. However, this notion can be introduced
and studied
for arbitrary nilpotent pairs (see \cite[sect.\,3]{pCM}, \cite{IMRN}): \par
$\bullet$ a nilpotent pair $\eb$ is called {\it rectangular\/}, if
$e_1$ and $e_2$ can be included in commuting $\tri$-triples.
\\[.5ex]
The following criterion was proved in \cite[1.9]{IMRN}:
\\[.5ex]
{\bf Proposition.} {\sl Let $\hb$ be a characteristic of a nilpotent pair
$\eb$. Then $\eb$ is rectangular if and only if $h_1\in\Ima(\ad e_1)$
if and only if $h_2\in\Ima(\ad e_2)$.}
\\[1ex]
Using this and our classification, it is easy to verify the following
claim.
\begin{s}{Theorem} \label{rect-distin}
Let $\eb$ be a distinguished pair in a classical simple Lie algebra $\g$
and $\Gamma$ the corresponding skew-graph. Then $\eb$ is rectangular if
and only if each connected component of $\Gamma$ is rectangular.
\end{s}%
This is a marvellous justification of our terminology\,!
\begin{s}{Corollary} \\
1. Let $\g$ be of type $\GR{B}{n}$ or $\GR{C}{n}$. Then each principal
pair in $\g$ is rectangular. \\
2. Let $\g$ be of type $\GR{D}{n}$. Then the non-rectangular principal
pairs correspond precisely to the near-rectangular skew-graphs.
In particular, if $\g$ has a non-rectangular principal pair, then
$n$ is odd.
\end{s}

\vno{3} \indent
{\footnotesize
{\bf A.E.}:\ \parbox[t]{160pt}{%
 {\it Razmadze Math. Institute \\
M.\,Aleksidze str., 1\\
380093 Tbilisi \\
Republic of Georgia}  \\
alela@rmi.acnet.ge}
\hspace{.5cm}
{\bf D.P.}:\ \parbox[t]{190pt}{%
 {\it Math. Department \\
M.I.R.E.A. \\
prosp. Vernadskogo, 78 \\
Moscow 117454 \quad Russia} \\ dmitri@panyushev.mccme.ru \\
panyush@dpa.msk.ru }
}
\end{document}